\input ./style/arxiv-general.cfg
\documentclass[MSNbibl,number,citesort,seceqn,dvips]{arxbj}
\makeatletter
   \@ifpackageloaded{graphicx}{}{\usepackage{graphicx}}
\makeatother

\volume{22}
\issue{2}
\pubyear{2016}
\firstpage{1055}
\lastpage{1092}
\doi{10.3150/14-BEJ686}
\docsubty{FLA}

\makeatletter
\newcommand{\ee}{\mathrm{e}\mathrm{e}}
\def\inplus{+}
\def\e{\mathrm{e}}
\newtheorem{theorem}{Theorem}[section]
\newtheorem{lemma}{Lemma}[section]
\newtheorem{proposition}[theorem]{Proposition}
\newtheorem{corollary}[theorem]{Corollary}
\newproclaim{definition}[theorem]{Definition}
\newremark{example}[theorem]{Example}
\newremark{remark}{Remark}
\makeatother

\begin{document}
\begin{frontmatter}

\title{Non-Gaussian semi-stable laws arising in sampling of finite point processes}
\runtitle{Semi-stable laws}

\begin{aug}
\author[A]{\inits{R.}\fnms{Ritwik}~\snm{Chaudhuri}\thanksref{e1}\ead[label=e1,mark]{ritwik@unc.edu}}
\and
\author[A]{\inits{V.}\fnms{Vladas}~\snm{Pipiras}\corref{}\thanksref{e2}\ead[label=e2,mark]{pipiras@email.unc.edu}}
\address[A]{Department of Statistics and OR, UNC-CH, Chapel Hill, NC
27599, USA.\\ \printead{e1,e2}}
\end{aug}

%
\received{\smonth{12} \syear{2013}}
%
\revised{\smonth{10} \syear{2014}}

%
\begin{abstract}
A finite point process is characterized by the distribution of the
number of points (the size) of the process. In some applications, for
example, in the context of packet flows in modern communication
networks, it is of interest to infer this size distribution from the
observed sizes of sampled point processes, that is, processes obtained
by sampling independently the points of i.i.d. realizations of the
original point process. A standard nonparametric estimator of the size
distribution has already been suggested in the literature, and has been
shown to be asymptotically normal under suitable but restrictive
assumptions. When these assumptions are not satisfied, it is shown here
that the estimator can be attracted to a semi-stable law. The
assumptions are discussed in the case of several concrete examples. A
major theoretical contribution of this work are new and quite general
sufficient conditions for a sequence of i.i.d. random variables to be
attracted to a semi-stable law.
\end{abstract}

%
\end{frontmatter}

\section{Introduction}
\label{sintro}
We first explain the motivation behind this work, namely, understanding
statistical properties of certain estimators arising when sampling
finite point process. The issues raised in the motivation require
developing new theoretical results on the domain of attraction of the
so-called semi-stable laws. We conclude this section by describing this
theoretical contribution, along with the structure of this work.

Let $W, W^{(i)}$, $i = 1, 2,\ldots, N$, be i.i.d. integer-valued
random variables with the probability mass function (p.m.f.)
$f_{W}(w)$, $w \geq1$. Let also $\operatorname{Bin}(n, q)$ denote a binomial
distribution with parameters $n\geq1$, $q\in(0,1)$. Consider random
variables $W_{q}, W_{q}^{(i)}$, $i=1, 2,\ldots, N$, obtained from $W,
W^{(i)}$, $i=1, 2,\ldots,N$, through the relationships $W_{q} =
\operatorname{Bin}(W, q)$ and $W_{q}^{(i)} = \operatorname{Bin}(W^{(i)}, q)$,
$i=1, 2,\ldots,N$
(independently across $i$). Note that $W_{q}$ takes values in $0, 1,
2,\ldots,W$. Let the probability mass function of $W_{q}$ be
$f_{W_{q}}(s)$, $s \geq0$. The basic interpretation of $W_{q}$ is as
follows. If an object consists of $W$ points (a finite point process)
and each point is sampled with a probability $q$, then the number of
sampled points is $W_{q}=\operatorname{Bin}(W, q)$.

One application of the above setting arises in modern communication
networks. A finite point process (an object) is associated with the
so-called packet flow (and a point is associated with a single packet).
Sampling is used in order to reduce the amount of data being collected
and processed. One basic problem that has attracted much attention
recently is the inference of $f_{W}$ from the observed sampled data
$W_{q}^{(i)}$, $i=1, 2,\ldots, N$ (in principle, $W_{q}^{(i)} = 0$ is
not observed directly, but the inference about the number of times
$W_{q}^{(i)} = 0$ is made through other means). See, for example,
Duffield, Lund and Thorup \cite{duffieldlundthorip2006}, Hohn and
Veitch \cite{hohnveitch2006}, Yang and Michailidis~\cite{yangmichailidis2007}.
For other, more recent progress on
sampling in communication networks, see Antunes and Pipiras \cite
{pipirasantunes2014a,pipirasantunes2014b}, and references therein.

We are interested here in some statistical properties of a
nonparametric estimator of $f_{W}(w)$, introduced in Hohn and Veitch
\cite{hohnveitch2006} and also considered in Antunes and Pipiras~\cite
{antunespipiras2011}. We first briefly outline how the estimator is
derived. Estimation of $f_{W}(w)$ is based on a theoretical inversion
of the relation
%
\begin{eqnarray}
\label{efwq-fw} f_{W_{q}}(s) &=& \sum_{w=s}^{\infty}
P(W_{q} = s|W=w)P(W = w)
\nonumber
\\[-8pt]
\\[-8pt]
\nonumber
& = &\sum_{w=s}^{\infty}
\pmatrix{w
\cr
s} q^{s}(1-q)^{w-s}f_{W}(w),\qquad s\geq0.
\end{eqnarray}
In terms of the moment generating functions $G_{W_{q}}(z) =
\sum_{s=0}^{\infty}z^{s}f_{W_{q}}(s)$ and $G_{W}(z) = \sum_{w=1}^{\infty}
z^{w}f_{W}(w)$, the relation
(\ref{efwq-fw}) can be written as $G_{W_{q}}(z) = G_{W}(zq + 1 - q)$.
By changing the variables $zq+ 1 - q = x$, one has $G_{W}(x) =
G_{W_{q}}(q^{-1}x - q^{-1}(1 - q))$ which has the earlier form but with
$q$ replaced by $q^{-1}$ (and $z$ replaced by $x$). This suggests that
(\ref{efwq-fw}) can be inverted as
%
\begin{eqnarray}
\label{efwq-fw-inversion} f_{W}(w) &=& \sum_{s=w}^{\infty}
\pmatrix{s
\cr
w}\bigl(q^{-1}\bigr)^{w}\bigl(1 - q^{-1}
\bigr)^{s-w}f_{W_{q}}(s)
\nonumber
\\[-8pt]
\\[-8pt]
\nonumber
&=&\sum_{s=w}^{\infty}
\pmatrix{s
\cr
w} \frac
{(-1)^{s-w}}{q^{s}}(1-q)^{s-w}f_{W_{q}}(s),\qquad w\geq1.
\end{eqnarray}
Antunes and Pipiras \cite{antunespipiras2011}, Proposition $4.1$,
showed that the inversion relation (\ref{efwq-fw-inversion}) holds when
%
\begin{equation}
\label{econd1} \sum_{s=n}^{\infty} \pmatrix{s
\cr
n} \frac{(1-q)^{s-n}}{q^{s}}f_{W_{q}}(s) = \sum_{w=n}^{\infty}
\pmatrix{w
\cr
n}2^{w-n}(1-q)^{w-n}f_{W}(w)<\infty,\qquad
n\geq1.
\end{equation}
Observe that (\ref{econd1}) always holds when $q\in(0.5,1)$. But when
$q\in(0, 0.5]$, the finiteness of the above expression depends on the
behavior of $f_{W}(w)$ as $w\rightarrow\infty$. We shall make the
assumption~(\ref{econd1}) throughout this work.

In view of (\ref{efwq-fw-inversion}), a natural nonparametric
estimator of $f_{W}$ is
%
\begin{equation}
\label{efWestimator} \widehat{f}_{W}(w) = \sum_{s=w}^{\infty}
\pmatrix{s
\cr
w}\frac
{(-1)^{s-w}}{q^{s}}(1-q)^{s-w}\widehat{f}_{W_{q}}(s),\qquad
w\geq1,
\end{equation}
where
%
\begin{equation}
\label{efWqempericalestimator} \widehat{f}_{W_{q}}(s) = \frac{1}{N}\sum
_{i=1}^{N} 1_{\{W_{q}^{(i)} =
s\}},\qquad s\geq0,
\end{equation}
is the empirical p.m.f. of $f_{W_{q}}$, and $1_{A}$ denotes the
indicator function of an event $A$. Note that, by using (\ref{efWestimator}) and (\ref{efwq-fw-inversion}),
%
\begin{equation}
\label{efhatw-fw} \sqrt{N}\bigl(\widehat{f}_{W}(w) - f_{W}(w)
\bigr) = \sum_{s=w}^{\infty}\pmatrix{s
\cr
w}
\frac{(-1)^{s-w}}{q^{s}}(1-q)^{s-w}\sqrt{N}\bigl(\widehat{f}_{W_{q}}(s) -
f_{W_{q}}(s)\bigr).
\end{equation}
Since
%
\begin{equation}
\label{efhatWq-fWq} \bigl\{\sqrt{N}\bigl(\widehat{f}_{W_{q}}(s) -
f_{W_{q}}(s)\bigr)\bigr\}_{s=0}^{\infty}
\mathop{\rightarrow}^{d} \bigl\{\xi(s)\bigr\}_{s=0}^{\infty},
\end{equation}
where $\{\xi(s)\}_{s=0}^{\infty}$ is a Gaussian process with zero mean
and covariance structure
\begin{eqnarray*}
E\bigl(\xi(s_{1})\xi(s_{2})\bigr) = f_{W_{q}}(s_{1})1_{\{s_{1} = s_{2}\}}
- f_{W_{q}}(s_{1})f_{W_{q}}(s_{2}),
\end{eqnarray*}
one may naturally expect that under suitable assumptions, (\ref
{efhatw-fw}) is asymptotically normal in the sense that
%
\begin{equation}
\label{efwhatconvergence} \bigl\{\sqrt{N}\bigl(\widehat{f}_{W}(w) -
f_{W}(w)\bigr)\bigr\}_{w=1}^{\infty}
\mathop{\rightarrow}^{d}\bigl\{S(\xi)_{w}\bigr\}_{w=1}^{\infty},
\end{equation}
where $\{S(\xi)_{w}\}_{w=1}^{\infty}$ is a Gaussian process. Antunes
and Pipiras \cite{antunespipiras2011}, Theorem $4.1$, showed
that~(\ref{efwhatconvergence}) holds indeed if $R_{q,w} <\infty, w\geq1$, where
%
\begin{eqnarray}
\label{eRqw} R_{q,w} &=&\sum_{s=w}^{\infty}
\pmatrix{s
\cr
w}^{2}\frac
{(1-q)^{2(s-w)}}{q^{2s}}f_{W_{q}}(s)
\nonumber
\\[-8pt]
\\[-8pt]
\nonumber
&=& \sum_{i=w}^{\infty} f_{W}(i)
(1-q)^{i-2w}\pmatrix{i
\cr
w}\sum_{s=w}^{i}
\pmatrix{s
\cr
w} \pmatrix{i-w
\cr
s-w}\bigl(q^{-1} - 1
\bigr)^{s}.
\end{eqnarray}
The quantity $R_{q,w}$ is naturally related to the limiting variance of
$\sqrt{N}\widehat f_{W}(w)$. Indeed, since $NE(\widehat
{f}_{W_{q}}(s_{1}) - f_{W_{q}}(s_{1}))(\widehat{f}_{W_{q}}(s_{2}) -
f_{W_{q}}(s_{2})) = f_{W_{q}}(s_{1})1_{\{s_{1}=s_{2}\}} -
f_{W_{q}}(s_{1})f_{W_{q}}(s_{2})$\vspace*{1pt} and by using (\ref{efhatw-fw})
and (\ref{efwq-fw-inversion}), the asymptotic variance of $\sqrt
{N}\widehat{f}_{W}(w)$ is expected to be $R_{q,w} - (f_{W}(w))^{2}$.
Requiring $R_{q,w}<\infty$ is then a natural assumption in proving (\ref
{efwhatconvergence}).

We are interested in $\widehat{f}_{W}(w)$ when the condition $R_{q,w}
<\infty, w\geq1$, is not satisfied. In fact, such a situation
is expected with many distributions. For example, we show in Section~\ref{smain-results-sampling} below that if $f_{W}(w) = (1-c)c^{w-1}$,
$w\geq1$, is a geometric distribution with parameter $c\in(0,1)$, then
the distribution of $f_{W_{q}}(s)$ is given by
%
\begin{eqnarray}
\label{egeometric} f_{W_{q}}(s) = %
\cases{ \displaystyle\frac{(1-q)(1-c)}{1 - c(1-q)},&\quad$\mbox{if } s=0$, \vspace*{2pt}
\cr
\displaystyle\frac{1}{c}c_{q}^{s}(1 -
c_{q}),&\quad$\mbox{if } s\geq1$,} %
\end{eqnarray}
where $c_{q} = \frac{cq}{1 - c(1-q)}$. Moreover, the condition
$R_{q,w}<\infty$ holds if and only if $c < \frac{q}{1-q}$ (see Section~\ref{smain-results-sampling}). Thus, for example, we are interested
what happens with $\widehat{f}_{W}(w)$ when $W_{q}$ has p.m.f. given by
(\ref{egeometric}) with $c \geq\frac{q}{1-q}$.

To understand what happens when $R_{q,w} = \infty$, observe from (\ref
{efWestimator}) and (\ref{efWqempericalestimator}) that
$\widehat{f}_{W}(w)$ can also be written as
%
\begin{equation}
\label{eempiricalfW} \widehat{f}_{W}(w) = \frac{1}{N}\sum
_{i=1}^{N}X_{i},
\end{equation}
where $X_{i}$, $i=1,2,\ldots,N$, are i.i.d. random variables defined as
%
\begin{equation}
\label{edescriptionxi} X_{i} = \pmatrix{W_{q}^{(i)}
\cr
w}
\frac{(-1)^{W_{q}^{(i)} -
w}}{q^{W_{q}^{(i)}}}(1-q)^{W_{q}^{(i)} - w}1_{\{W_{q}^{(i)} \geq w\}}.
\end{equation}
Focus on the key term $\frac{(1 -
q)^{W_{q}^{(i)}}}{q^{W_{q}^{(i)}}}=(q^{-1} - 1)^{W_{q}^{(i)}}$ entering
(\ref{edescriptionxi}). For example, when $W$ is geometric with
parameter $c$, $W_{q}^{(i)}$ has p.m.f. in (\ref{egeometric}). One
then expects that
%
\begin{eqnarray}
\label{eheavytail} P\bigl(\bigl(q^{-1} -1 \bigr)^{W_{q}^{(i)}} > x
\bigr)&=&P \biggl(W_{q}^{(i)} > \frac{\operatorname
{log } x}{\operatorname{log }(q^{-1} - 1)} \biggr)
\nonumber
\\[-8pt]
\\[-8pt]
\nonumber
&\approx& \frac{1}{c}c_{q}^{{\operatorname{log } x}/{\operatorname{log }(q^{-1} -
1)}}=\frac{1}{c}
x^{-\alpha},
\end{eqnarray}
where $\alpha= \frac{\operatorname{log }c_{q}^{-1}}{\operatorname{log }(q^{-1} - 1)}$.
This suggests that the distribution of $X_{i}$, $i=1, 2,\ldots,N$, has
heavy tail and that the estimator $\widehat{f}_{W}(w)$ is
asymptotically non-Gaussian stable when $\alpha< 2$. In fact, the
story turns out to be more complex. Because of the discrete nature of
$W_{q}^{(i)}$, the relation (\ref{eheavytail}) does not hold in the
asymptotic sense as $x\rightarrow\infty$. An appropriate setting in
this case involves the so-called semi-stable laws. In the semi-stable
context, moreover, the convergence of (\ref{eempiricalfW}) is
expected only along subsequences of $N$.

Semi-stable laws have been studied quite extensively (see Section
\ref{ssemistable} for references). They are infinitely divisible and extend the
stable laws by allowing the power function in the L\'evy measure (of
the stable law) to be multiplied by a function with a multiplicative
period. In particular, necessary and sufficient conditions are known
for a distribution to be attracted to a semi-stable law (see Theorem
\ref{tmain-semistable} below), that is, for the sum of independent copies following the
distribution to converge to a semi-stable law (along a subsequence and
after suitable normalization and centering). A common example (and, in
fact, one of the few concrete examples) of such a distribution is that
of a log-geometric random variable
%
\begin{equation}
\label{enewadded} X = a^{W_{q}} \qquad\mbox{with } P(W_{q} = s) =
(1-c_{q})c_{q}^{s}, s = 0,1,\ldots,
\end{equation}
where $a>0$ and $c_{q}\in(0,1)$. (Strictly speaking, the log-geometric
case is when $a=e$.) Note that in (\ref{enewadded}), we use purposely
the notation of (\ref{egeometric}) and (\ref{edescriptionxi}).

In fact, motivated by (\ref{edescriptionxi}) and the desire to
consider more general distributions than log-geometric, we will show
that the domain of attraction of semi-stable laws also includes the
distributions of random variables of the form
%
\begin{equation}
\label{enewdded2} X = k(W_{q})a^{W_{q}} \qquad\mbox{with }
P(W_{q} = s) = h(s)c_{q}^{s}, s = 0,1,\ldots,
\end{equation}
where $k$ and $h$ are functions satisfying suitable but also flexible
conditions. Our approach goes through verifying that the distributions
determined by (\ref{enewdded2}) satisfy the necessary and sufficient
conditions to be attracted to a semi-stable law. Somewhat surprising
perhaps, the proof turns out to be highly nontrivial. The difficulty
lies in dealing with the general case when both functions $k$ and $h$
in (\ref{enewdded2}) are not constant. Much of this work, in fact,
concerns this problem.

The rest of this work is structured as follows. Preliminaries on
semi-stable laws can be found in Section~\ref{ssemistable}. In Section~\ref{smain-results}, we state and prove the main general results of
this work concerning semi-stable distributions and their domains of
attraction. In Section~\ref{smain-results-sampling}, we apply the main
results from Section~\ref{smain-results} to sampling of finite point
processes. Several concrete examples, in particular, are considered. A
few auxiliary results are given in the \hyperref[sauxiliary]{Appendix}. Some
numerical illustrations can be found in Chaudhuri and Pipiras \cite
{chaudhuripipiras2013}.

\section{Preliminaries on semi-stable laws}
\label{ssemistable}
One way to characterize a semi-stable distribution is through its
characteristic function (Maejima~\cite{maejima2001}).

\begin{definition}
A probability distribution $\mu$ on $\mathbb{R}$ (or a random variable
with distribution $\mu$) is called semi-stable if there exist $r, b \in
(0,1)$ and $c \in\mathbb{R}$ such that
%
\begin{equation}
\label{echaracteristicsemistable} \widehat{\mu}(\theta)^{r} = \widehat{\mu}(b
\theta)\e^{\mathrm{i}c\theta}\qquad\mbox{for all } \theta\in\mathbb{R},
\end{equation}
and $\widehat{\mu}(\theta) \neq0,\mbox{ for all }\theta\in\mathbb{R}$,
where $\widehat{\mu}(\theta)$ denotes the characteristic function of
$\mu$.
\end{definition}

A semi-stable distribution is known to be infinitely divisible
(Maejima \cite{maejima2001}) with a location parameter $\eta\in\mathbb
{R}$, a Gaussian part with variance $\sigma^{2}\geq0$ and a
non-Gaussian part with L\' evy measure characterized by (distribution) functions
%
\begin{equation}
\label{elevyleft} L(x) = \frac{M_{L}(x)}{|x|^{\alpha}},\qquad x<0,\qquad
 R(x) = -\frac{M_{R}(x)}{x^{\alpha}},\qquad  x>0,
\end{equation}
where $\alpha\in(0,2)$, $M_{L}(c^{{1}/{\alpha}} x) = M_{L}(x)$ when
$x < 0$, and $M_{R}(c^{{1}/{\alpha}} x) = M_{R}(x)$ when $x > 0$,
for some $c > 0$. The functions $M_{L}$ and $M_{R}$ are thus periodic
with multiplicative period $c^{{1}/{\alpha}}$. The functions $L(x)$
and $R(x)$ are left-continuous and non-decreasing on $(-\infty, 0)$ and
right-continuous and non-decreasing on $(0, \infty)$, respectively.
The characteristic function of a semi-stable distribution with a
location parameter $\eta$ and without a Gaussian part is given by
%
\begin{equation}
\label{echaracteristic} \log\widehat{\mu}(t) = \mathrm{i}\eta t + \int_{-\infty}^{0}
\biggl(\e^{\mathrm{i}tx} - 1 - \frac{\mathrm{i}tx}{1+x^{2}}\biggr)\,\mathrm{d}L(x) + \int
_{0}^{\infty} \biggl(\e^{\mathrm{i}tx} - 1 -
\frac
{\mathrm{i}tx}{1+x^{2}}\biggr)\,\mathrm{d}R(x).
\end{equation}

Semi-stable distributions arise as limits of partial sums of
i.i.d. random variables. Let $X_{1}, X_{2},\ldots$ be a sequence of
i.i.d. random variables with a common distribution function $F$.
Consider the sequence of partial sums
%
\begin{equation}
\label{enormalizedsum} S^{*}_{n} = \frac{1}{A_{k_{n}}} \Biggl\{\sum
_{j=1}^{k_{n}} X_{j} -
B_{k_{n}} \Biggr\},
\end{equation}
where $\{A_{k_{n}}\}$ and $\{B_{k_{n}}\}$ are normalizing and centering
sequences. Semi-stable laws arise as limits of partial sums
$S^{*}_{n}$, supposing that $\{k_{n}\}$ satisfies
%
\begin{equation}
\label{econdkn} k_{n}\rightarrow\infty, k_{n}\leq
k_{n+1},\qquad \lim_{n\rightarrow\infty}\frac{k_{n+1}}{k_{n}}= c\in[1,
\infty).
\end{equation}
Moreover, if $S_{n}^{*}$ converges to a nontrivial limit (semi-stable
distribution), the distribution $F$ of $X_{j}$ is said to be in the
domain of attraction of the limiting semi-stable law. In this case and
supposing the limiting law is non-Gaussian semi-stable, it is known
that the normalizing sequence $\{A_{k_{n}}\}$ necessarily satisfies
%
\begin{equation}
\label{econditionAkn} A_{k_{n}}\rightarrow\infty, A_{k_{n}}\leq
A_{k_{n+1}}, \qquad\lim_{n\rightarrow\infty}\frac{A_{k_{n+1}}}{A_{k_{n}}}=
c^{{1}/{\alpha}} \qquad\mbox{where }\alpha\in(0, 2).
\end{equation}

Megyesi \cite{megyesi2000}, Grinevich and Khokhlov \cite
{grinevichkhokhlov1995} gave necessary and sufficient conditions for
a distribution to be in the domain of attraction of a semi-stable distribution.

\begin{theorem}[(Megyesi \cite
{megyesi2000}, Corollary~3)]\label{tmain-semistable}
Distribution $F$ is in the domain of attraction of a non-Gaussian
semi-stable distribution with the characteristic function (\ref
{echaracteristic}) along the subsequence $k_{n}$ with normalizing
constants ${A_{k_{n}}}$ satisfying (\ref{econdkn}) and (\ref
{econditionAkn}) if and only if for all $x > 0$ large enough,
%
\begin{eqnarray}
\label{ecorollary3,megyesi,lefttail} F_{- }(-x) &=& x^{-\alpha}l^{*}(x)
\bigl(M_{L}\bigl(-\delta(x)\bigr) + h_{L}(x)\bigr),
\\
\label{ecorollary3,megyesi,righttail} 1 - F(x) &=& x^{-\alpha}l^{*}(x)
\bigl(M_{R}\bigl(\delta(x)\bigr) + h_{R}(x)\bigr),
\end{eqnarray}
where $l^{*}$ is a right-continuous function, slowly varying at $\infty
$, $\alpha\in(0,2)$, $F_{-}$ is the left-continuous version of $F$ and
the error functions $h_{R}$ and $h_{L}$ are such that
%
\begin{equation}
\label{eerrorfunction} h_{K}(A_{k_{n}}x_{0})\rightarrow0\qquad
\mbox{as } n\rightarrow\infty,
\end{equation}
for every continuity point $x_{0}$ of $M_{R}$, if $K = R$, and $-x_{0}$
of $M_{L}$, if $K = L$. $M_{K}$, $K\in\{L,R\}$, are two periodic
functions with common multiplicative period $c^{{1}/{\alpha}}$ and
for all large enough $x$, $\delta(x)$ is defined as
%
\begin{equation}
\label{econditiondeltax} \delta(x) = \frac{x}{a(x)}\in\bigl[1, c^{{1}/{\alpha}} +
\varepsilon\bigr],
\end{equation}
where $\varepsilon> 0$ is any fixed number, with
%
\begin{equation}
\label{econditionax} a(x) = A_{k_{n}}\qquad \mbox{if } A_{k_{n}} \leq x <
A_{k_{n+1}}.
\end{equation}
\end{theorem}

Grinevich and Khokhlov \cite{grinevichkhokhlov1995} also showed that,
in the sufficiency part of the theorem above, $k_{n}$ can be chosen as
follows. First, choose a sequence $\{\tilde{A}_{n}\}$ such that
%
\begin{equation}
\label{edefineAntilde} \lim_{n\rightarrow\infty}n\tilde{A}_{n}^{-\alpha}l^{*}(
\tilde{A}_{n}) = 1
\end{equation}
and
%
\begin{equation}
\label{econditionAntilde} \tilde{A}_{n}\rightarrow\infty,\qquad \tilde{A}_{n}
\leq\tilde {A}_{n+1}\quad \mbox{and}\quad \lim_{n\rightarrow\infty}
\frac{\tilde
{A}_{n+1}}{\tilde{A}_{n}}=1.
\end{equation}
Define a new sequence $\{a_{n}\}$ by setting $a_{n} = A_{k_{n}}$ for
every $n$, where $A_{k_{n}}$ appears in (\ref{econditionax}). Then,
the natural numbers $k_{n}$ can be chosen as
%
\begin{equation}
\label{edefinitionAk{n}} \tilde{A}_{k_{n}}\leq a_{n} <
\tilde{A}_{k_{n+1}}.
\end{equation}
The centering constants $B_{k_{n}}$ in (\ref{enormalizedsum}) can be
chosen as (Cs\"{o}rg\"{o} and Megyesi \cite{csorgomegyesi2003})
%
\begin{equation}
\label{edefiBkn} B_{k_{n}} = k_{n}\int_{{1}/{k_{n}}}^{1 - {1}/{k_{n}}}
Q(s) \,\mathrm{d}s,
\end{equation}
where, for $0\leq s\leq1$,
%
\begin{equation}
\label{edefinitionQs} Q(s) = \inf_{y}\bigl\{F(y)\geq s\bigr\}.
\end{equation}
The location parameter $\eta$ of the limiting semi-stable law in (\ref
{echaracteristic}) is then given by
%
\begin{eqnarray}
\label{elocationoriginal} \eta= \Theta(\psi_{1}) - \Theta(\psi_{2}),
\end{eqnarray}
where
%
\begin{equation}
\label{eaddedcharcteristic} \Theta(\psi_{i}) = \int_{0}^{1}
\frac{\psi_{i}(s)}{1 + \psi_{i}^{2}(s)} \,\mathrm{d}s - \int_{1}^{\infty}
\frac{\psi_{i}^{3}(s)}{1 + \psi_{i}^{2}(s)} \,\mathrm{d}s,\qquad i = 1,2,
\end{equation}
and
%
\begin{equation}
\label{elevyfunctioninversion} \psi_{1} (s) = \inf_{x < 0}\bigl\{L(x)
> s\bigr\}, \qquad\psi_{2}(s) = \inf_{x<0}\bigl\{-R(-x) >
s\bigr\}.
\end{equation}
It is also worth mentioning that the slowly varying function $l^{*}(x)$
entering in (\ref{ecorollary3,megyesi,lefttail}) and (\ref
{ecorollary3,megyesi,righttail}) can be replaced by two different,
asymptotically equivalent slowly varying functions $l_{1}^{*}(x)$ and
$l_{2}^{*}(x)$. The proof of this result is given in Lemma~\ref
{ltwo-slow-var} in the \hyperref[sauxiliary]{Appendix}.
\section{General results concerning semi-stable domain of attraction}
\label{smain-results}
The next theorem is the main result of this work. We use the following
notation throughout this work:
\begin{eqnarray*}
\lceil x \rceil&=& \mbox{ the smallest integer larger than or equal to } x,
\\
\lceil x \rceil_{\inplus} &=& \mbox{ the smallest integer strictly larger
than } x.
\end{eqnarray*}
For example, $\lceil2.47 \rceil= \lceil2.47 \rceil_{\inplus} = 3$
but $\lceil3 \rceil=3$ and $\lceil3 \rceil_{\inplus} = 4$. The
function $\lceil x \rceil_{\inplus}$ is the right-continuous version of
the function $\lceil x \rceil$. Also note that $\lceil x \rceil
_{\inplus} = [x]+1$, where $[x]$ is the integer part of $x$ (i.e.,
the largest integer smaller than or equal to $x$).

\begin{theorem}\label{tmain-th}
Let $W_{q}$ be an integer-valued random variable taking values in $0,
1, 2, \ldots$ such that, for all $x > 0$,
%
\begin{equation}
\label{eevenWq} P \biggl(\frac{W_q}{2} \geq x, W_{q}\mbox{ is
even} \biggr) = \sum_{n =
\lceil x\rceil}^{\infty} P \biggl(
\frac{W_q}{2} = n \biggr) = h_{1}\bigl(\lceil x \rceil
\bigr)\e^{-\nu\lceil x \rceil},
\end{equation}
%
\begin{equation}
\label{eoddWq} P \biggl(\frac{W_q - 1}{2} \geq x, W_{q} \mbox{ is
odd} \biggr) = \sum_{n
= \lceil x\rceil}^{\infty} P \biggl(
\frac{W_q - 1}{2} = n \biggr) = h_{2}\bigl(\lceil x \rceil
\bigr)\e^{-\nu\lceil x \rceil},
\end{equation}
where $\nu> 0$ and the functions $h_{1}$ and $h_{2}$ satisfy
%
\begin{equation}
\label{eh2h1ratio} \frac{h_{2}(x)}{h_{1}(x)} \rightarrow c_{1}\qquad \mbox{as } x
\rightarrow \infty,
\end{equation}
for some fixed $c_{1}\geq0$, and
%
\begin{equation}
\label{eh1ratio} \frac{h_{1}(ax)}{h_{1}(x)}\rightarrow1 \qquad\mbox{as }
x\rightarrow\infty, a
\rightarrow1.
\end{equation}
Let also
%
\begin{equation}
\label{erandomvariabledefinition} X = L\bigl(\e^{W_{q}}\bigr)\e^{\beta W_{q}}(-1)^{W_{q}},
\end{equation}
where $\beta> 0$ and $L$ is a slowly varying function at $\infty$ such
that $L(\e^{n})$ is ultimately monotonically increasing. Suppose that
%
\begin{equation}
\alpha:= \frac{\nu}{2\beta} < 2.
\end{equation}
Then, $X$ is attracted to the domain of a semi-stable distribution in
the following sense. If $X, X_{1}, X_{2},\ldots$ are i.i.d. random
variables, then
as $n\rightarrow\infty$, the partial sums
%
\begin{equation}
\label{epartialsum} \frac{1}{A_{k_{n}}} \Biggl\{\sum_{j=1}^{k_{n}}
X_{j} - B_{k_{n}} \Biggr\}
\end{equation}
converge to a semi-stable distribution with
%
\begin{equation}
\label{edefinitionknAknBkn} k_{n}= \biggl\lceil\frac{\e^{(n-1)\nu}}{h_{1}(n-1)}
\biggr\rceil,\qquad
A_{k_{n}}= L\bigl(\e^{2n-2}\bigr)\e^{2\beta(n-1)}
\end{equation}
and $B_{k_{n}}$ given by (\ref{edefiBkn}).
The limiting semi-stable distribution is non-Gaussian, has location
parameter given in (\ref{elocationoriginal}) and is characterized by
%
\begin{eqnarray}
\label{echarcterization0} \alpha= \frac{\nu}{2\beta},
\end{eqnarray}
%
\begin{eqnarray}
\label{echarcterization} M_{L}(-x) &=& c_{1}\e^{-\nu([{1}/{2}+({1}/{(2\beta)})\log x] -
({1}/{(2\beta)})\log x)},
\nonumber
\\[-8pt]
\\[-8pt]
\nonumber
M_{R}(x) &=& \e^{-\nu(\lceil({1}/{(2\beta)})\log x\rceil_{\inplus} -
({1}/{(2\beta)})\log x)},\qquad x > 0.
\end{eqnarray}
\end{theorem}

\begin{pf} The result will be proved by verifying the sufficient
conditions (\ref{ecorollary3,megyesi,lefttail})--(\ref{ecorollary3,megyesi,righttail}) of Theorem~\ref{tmain-semistable}. We break the
proof into two cases dealing with (\ref{ecorollary3,megyesi,lefttail}) and (\ref{ecorollary3,megyesi,righttail})
separately. The final part of the proof shows that the sequence $k_{n}$
can be chosen as in (\ref{edefinitionknAknBkn}).

\textit{Step} 1 (showing (\ref{ecorollary3,megyesi,righttail})): Fix $x>0$ large enough. In view of (\ref
{erandomvariabledefinition}), we are interested in
%
\begin{equation}
\label{eFxrighttail} \bar{F}(x) := 1 - F(x) = P \bigl(L \bigl(\e^{W_{q}}
\bigr)\e^{\beta
W_{q}}(-1)^{W_{q}} > x \bigr).
\end{equation}
Let $Z_{2} = \frac{W_{q}}{2}$. Note that (\ref{eFxrighttail}) can be
written as
%
\begin{eqnarray}
\label{ecompactrighttail} \bar{F}(x)&=& P \bigl(L\bigl(\e^{2Z_{2}}\bigr)\e^{2\beta Z_{2}}
> x, Z_{2}\mbox{ is integer} \bigr)
\nonumber
\\
&=&P \bigl(L\bigl(\e^{2Z_{2}}\bigr)\e^{2\beta Z_{2}} > x \bigr)
\\
&=&P \biggl(Z_{2} + \frac{1}{2\beta} \log L \bigl(\e^{2Z_{2}}
\bigr) > \frac
{1}{2\beta}\log x \biggr),\nonumber
\end{eqnarray}
where, in view of (\ref{eevenWq}),
%
\begin{equation}
\label{edefinitionz2} P(Z_{2} \geq x) = h_{1}\bigl(\lceil x \rceil
\bigr)\e^{-\nu\lceil x \rceil}.
\end{equation}
We next want to write $\bar{F}(x)$ in (\ref{ecompactrighttail}) as
%
\begin{equation}
\label{ecompactrighttaill} \bar{F}(x)=P \biggl(Z_{2} \geq g \biggl(
\frac{1}{2\beta}\log x \biggr) \biggr)
\end{equation}
for some function $g$.

There are many choices for $g$ in (\ref{ecompactrighttaill}). One
natural choice is to take
%
\begin{equation}
\label{edefinitiong0} g_{0}(y) = n\qquad \mbox{if } (n-1) + \frac{1}{2\beta}
\log L \bigl(\e^{2n-2} \bigr) \leq y < n + \frac{1}{2\beta}\log L
\bigl(\e^{2n} \bigr).
\end{equation}
The function $g_{0}$, however, turns out not to be suitable for our
purpose. It will be used below only for reference and comparison to
other related functions. We will use a related function $g_{1}$
defined, for integer $n \geq2$, as
%
\begin{equation}
\label{edefinitiong1} g_{1}(y) =
\cases{n - 1, \vspace*{2pt}\cr
\quad \mbox{if } \displaystyle n - 1 + \frac{1}{2\beta} \log L \bigl(\e^{2n-2} \bigr)
\leq y < n - 1 +
\frac{1}{2\beta} \log L \bigl(\e^{2n} \bigr),\vspace*{2pt}
\cr
\displaystyle y -
\frac{1}{2\beta} \log L \bigl(\e^{2n} \bigr),\vspace*{2pt}\cr
\quad \mbox{if }\displaystyle n-1 +
\frac
{1}{2\beta} \log L \bigl(\e^{2n} \bigr) \leq y < n +
\frac{1}{2\beta} \log L \bigl(\e^{2n} \bigr).} 
\end{equation}
We will also use the function
%
\begin{equation}
\label{edefinitiong2} g_{2}(y) = f^{-1}(y) = \inf\bigl\{z\dvt f(z)
\geq y\bigr\}
\end{equation}
defined as an inverse of the function
%
\begin{equation}
\label{edefinitionf} f(z) = z + \frac{1}{2\beta}\log L \bigl(\e^{2z} \bigr).
\end{equation}
Note that
%
\begin{eqnarray}
\label{eequality} \bigl\lceil g_{0}(y) \bigr\rceil= \bigl\lceil
g_{1}(y) \bigr\rceil_{\inplus} = \bigl\lceil g_{2}(y)
\bigr\rceil_{\inplus} = \bigl\lceil g(y) \bigr\rceil,
\end{eqnarray}
where $g$ is any function satisfying (\ref{ecompactrighttaill}). The
functions $g_{0}$, $g_{1}$ and $g_{2}$ are plotted in Figure~\ref{Figure1}.
%
\begin{figure}

\includegraphics{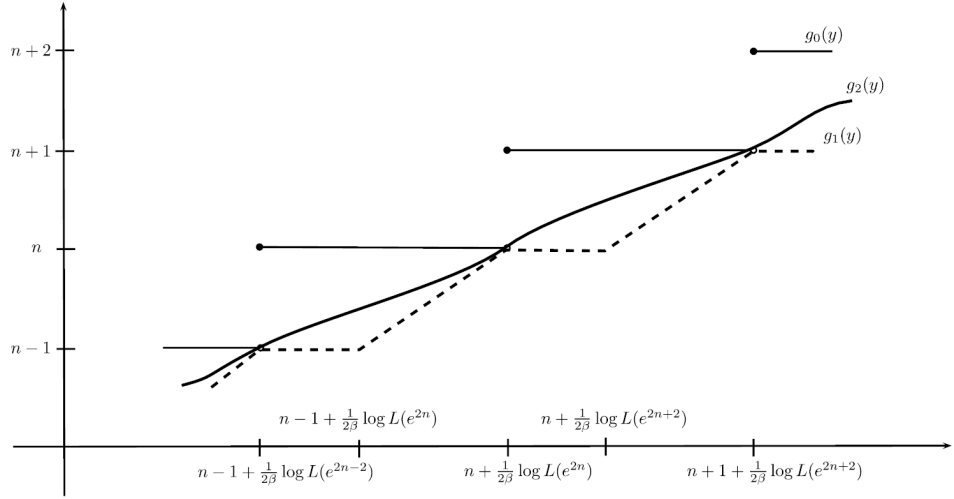}

\caption{Plot of $g_{0}(y)$, $g_{1}(y)$ and $g_{2}(y)$.}
\label{Figure1}
\end{figure}

We shall use another function $\tilde{g}_{1}$ which modifies
$g_{1}$ in the following way:
for $n \geq2$,
%
\begin{eqnarray}
\label{edefinitiong1tilde} &&\tilde{g}_{1}(y) = y - \frac{1}{2\beta}\log L
\bigl(\e^{2n-2}\bigr)
\nonumber
\\[-8pt]
\\[-8pt]
\eqntext{\mbox{if }\displaystyle n-1 + \frac{1}{2\beta}\log L
\bigl(\e^{2n-2}\bigr) \leq y < n + \frac{1}{2\beta
}\log L
\bigl(\e^{2n}\bigr).\qquad\qquad}
\end{eqnarray}
One relationship between the functions $g_{1}$ and $\tilde{g}_{1}$ can
be found in Lemma~\ref{lg1-tilde-g1} in the \hyperref[sauxiliary]{Appendix},
and will be used in the proof below.
Note that $\tilde{g}_{1}(y)$ can be expressed as
%
\begin{equation}
\label{eg1definedasdifference} \tilde{g}_{1}(y) = y - \tilde{g}_{1}^{*}(y),
\end{equation}
where, for $n \geq2$,
%
\begin{equation}
\label{etildeg1*} \tilde{g}_{1}^{*}(y) = \frac{1}{2\beta}
\log L \bigl(\e^{2n-2} \bigr)\qquad \mbox{if } n-1 + \frac{1}{2\beta}\log L
\bigl(\e^{2n-2}\bigr) \leq y < n + \frac
{1}{2\beta}\log L
\bigl(\e^{2n}\bigr).
\end{equation}
See Lemma~\ref{lg-tilde-*} in the \hyperref[sauxiliary]{Appendix} for a
property of $\tilde{g}_{1}^{*}$ which will be used in the proof below.\vadjust{\goodbreak}

We need few properties of the function $g_{2}$. Since $g_{2}$ is
the inverse of the function $f$, we have $\e^{g_{2} (\operatorname{log }
x )}$ as the inverse of $\e^{f (\operatorname{log } x )}$. Indeed,
\[
\e^{g_{2} (\operatorname{log }\e^{f (\operatorname{log } x )} )} = \e^{g_{2}(f(\operatorname{log }x))} = \e^{\operatorname{log } x} = x.
\]
Note now from (\ref{edefinitionf}) that
\[
\label{edefinitioneh*logx} \e^{f (\operatorname{log } x )} = \e^{\operatorname{log }x
+({1}/{(2\beta)})\operatorname
{log } L(x^{2})} = x \bigl(L
\bigl(x^{2}\bigr) \bigr)^{{1}/{(2\beta)}}.
\]
Since $ (L(x^2) )^{{1}/{(2\beta)}}$ is a slowly varying
function, $\e^{f (\operatorname{log } x )}$ is a regularly varying
function. So, by Theorem~1.5.13 of Bingham, Goldie and Teugels \cite
{binghamgoldieteugels1987},
\[
\label{eeg2logx} \e^{g_{2} (\operatorname{log } x )} = x l(x),
\]
where $l(x)$ is a slowly varying function. Hence,
\[
\label{eg2logx} g_{2} (\operatorname{log }x ) = \operatorname{log }x +
\operatorname{log }l(x) = \operatorname {log }x + g_{2}^{*}
(\operatorname{log } x ),
\]
where
\[
\label{g2*logx} g_{2}^{*} (\operatorname{log }x ) =
\operatorname{log }l(x)
\]
or replacing log $x$ by $y$,
%
\begin{equation}
\label{eg2representation} g_{2}(y) = y + g_{2}^{*}(y).
\end{equation}
Note also that for any $A>0$, we have
%
\begin{equation}
\label{eg2*property} g^{*}_{2} (\operatorname{log }Ax ) -
g^{*}_{2} (\operatorname{log }x ) = \operatorname{log
}l(Ax) - \operatorname{log }l(x)
= \operatorname{log }\frac{l(Ax)}{l(x)} \rightarrow0\qquad \mbox{as } x\rightarrow
\infty.
\end{equation}

Continuing with (\ref{ecompactrighttaill}) now, note that, by
using (\ref{edefinitionz2}) and (\ref{eequality}),
%
\begin{eqnarray}
\label{eexpressedrighttail} \bar{F}(x)&= &P \biggl(Z_{2} \geq g \biggl(
\frac{1}{2\beta}\log x \biggr) \biggr)
\nonumber
\\
&=& h_{1} \biggl(\biggl\lceil g\biggl(\frac{1}{2\beta}\log x\biggr)
\biggr\rceil \biggr)\e^{-\nu
\lceil g(({1}/{(2\beta)})\log x)\rceil}
\\
&=& h_{1} \biggl(\biggl\lceil g_{2}\biggl(
\frac{1}{2\beta}\log x\biggr)\biggr\rceil_{\inplus} \biggr)
\e^{-\nu\lceil g_{1}(({1}/{(2\beta)})\log x)\rceil_{\inplus}}.\nonumber
\end{eqnarray}
By using (\ref{eg1definedasdifference}), note further that
%
\begin{eqnarray}
\label{eexpressedrighttaill} \bar{F}(x) &=& h_{1} \biggl(\biggl\lceil
g_{2}\biggl(\frac{1}{2\beta}\log x\biggr)\biggr\rceil
_{\inplus} \biggr)\e^{-\nu\tilde{g}_{1}(({1}/{(2\beta)})\log x)}
\nonumber
\\
&& {}\times\e^{-\nu({g}_{1}(({1}/{(2\beta)})\log
x) - \tilde{g}_{1}(({1}/{(2\beta)})\log x))}\e^{-\nu(\lceil g_{1}(
({1}/{(2\beta)})\log x)\rceil_{\inplus} - g_{1}(({1}/{(2\beta)})\log
x))}\nonumber
\\
&=&h_{1} \biggl(\biggl\lceil g_{2}\biggl(
\frac{1}{2\beta}\log x\biggr)\biggr\rceil_{\inplus} \biggr)
\e^{-\nu(({1}/{(2\beta)})\log x - \tilde{g}_{1}^{*}(({1}/{(2\beta)})\log x))}
\\
&& {}\times\e^{-\nu({g}_{1}(({1}/{(2\beta)})\log
x) - \tilde{g}_{1}(({1}/{(2\beta)})\log x))}\e^{-\nu(\lceil g_{1}(
({1}/{(2\beta)})\log x)\rceil_{\inplus} - g_{1}(({1}/{(2\beta)})\log
x))}
\nonumber
\\
&=&x^{-\alpha}l_{1}^{*}(x) \bigl(M_{R}
\bigl(\delta(x)\bigr) + h_{R}(x)\bigr),\nonumber
\end{eqnarray}
where $\alpha= \frac{\nu}{2\beta}$ as given in (\ref{echarcterization0}),
%
\begin{eqnarray}
\label{el1*x} l_{1}^{*}(x) &=& h_{1} \biggl(
\biggl\lceil g_{2}\biggl(\frac{1}{2\beta}\log x\biggr)\biggr\rceil
_{\inplus} \biggr)\e^{\nu\tilde{g}_{1}^{*}(({1}/{(2\beta)})\log x)}
\nonumber
\\[-8pt]
\\[-8pt]
\nonumber
&&{}\times\e^{-\nu
({g}_{1}(({1}/{(2\beta)})\log x) - \tilde{g}_{1}(({1}/{(2\beta)})\log x))},
\\
\label{edefinitionMRdeltax} M_{R}\bigl(\delta(x)\bigr) &=& \e^{-\nu(\lceil\tilde{g}_{1}(({1}/{(2\beta)})\operatorname
{log }x)\rceil_{\inplus} - \tilde{g}_{1}(({1}/{(2\beta)})\operatorname{log }x))}
\end{eqnarray}
and
%
\begin{eqnarray}
\label{edefinitionhR} h_{R}(x) &=& \e^{-\nu(\lceil g_{1}(({1}/{(2\beta)})\operatorname{log }x)\rceil
_{\inplus} - g_{1}(({1}/{(2\beta)})\operatorname{log }x))}
\nonumber
\\[-8pt]
\\[-8pt]
\nonumber
&&{}- \e^{-\nu(\lceil
\tilde{g}_{1}(({1}/{(2\beta)})\operatorname{log }x)\rceil_{\inplus} - \tilde
{g}_{1}(({1}/{(2\beta)})\operatorname{log }x))}.
\end{eqnarray}
We next show that the functions $l_{1}^{*}$, $M_{R}$ and $h_{R}$
satisfy the conditions of Theorem~\ref{tmain-semistable} with suitable
choices of $\delta(x)$ and $A_{k_{n}}$.

By Lemma~\ref{ll1*} in the \hyperref[sauxiliary]{Appendix}, $l_{1}^{*}(x)$
is a right-continuous slowly varying function and hence it satisfies
the conditions of Theorem~\ref{tmain-semistable}. For the function
$M_{R}(\delta(x))$, note from (\ref{edefinitionMRdeltax}) that
%
\begin{eqnarray}
\label{erecursiveformMR} M_{R}\bigl(\delta(x)\bigr) &=&\e^{-\nu (
\lceil{2\beta\tilde
{g}_{1}(({1}/{(2\beta)})\operatorname{log }x)}/{(2\beta)} \rceil_\inplus -
{2\beta\tilde{g}_{1}(({1}/{(2\beta)})\operatorname{log }x)}/{(2\beta)}
)}
\nonumber
\\[-8pt]
\\[-8pt]
\nonumber
&=&M_{R} \bigl(\e^{2\beta\tilde{g}_{1}(({1}/{(2\beta)})\operatorname{log }x)} \bigr)
\end{eqnarray}
with
%
\begin{equation}
\label{edefinitionMR} M_{R}(x) =
 \e^{-\nu(\lceil{\operatorname{log } x}/{(2\beta)}\rceil_{\inplus} -
{\operatorname{log }x}/{(2\beta)})}.
\end{equation}
The function $M_{R}(x)$ is periodic with multiplicative period
$\e^{2\beta}$, and is right-continuous as required in Theorem~\ref
{tmain-semistable}. Since the period $\e^{2\beta}$ is also $c^{{1}/{\alpha}}$, this yields
%
\begin{equation}
\label{edefinitionc} c = \e^{\nu}.
\end{equation}
To choose $\delta(x)$, note from (\ref{erecursiveformMR}) that
\[
M_{R}\bigl(\delta(x)\bigr) =M_{R} \bigl(\e^{2\beta\tilde{g}_{1}(({1}/{(2\beta)})\operatorname{log }x) - 2\beta(n-1)}
\bigr),
\]
for any $n\geq1$, since $M_{R}$ has multiplicative period $\e^{2\beta
}$. We can set
%
\begin{equation}
\label{edefinitiondeltax} \delta(x) = \e^{2\beta\tilde{g}_{1}(({1}/{(2\beta)})\operatorname{log }x) -
2\beta(n-1)}\qquad \mbox{if }\e^{2\beta(n-1)}L
\bigl(\e^{2n-2}\bigr)\leq x <\e^{2n\beta}L\bigl(\e^{2n}\bigr).
\end{equation}
From (\ref{edefinitiong1tilde}), we have
%
\begin{eqnarray}
\label{edeltaxfinalform} \delta(x)&=& \e^{2\beta(({1}/{(2\beta)})\operatorname{log }
 x -( {1}/{(2\beta)})\operatorname{log }L(\e^{2n-2})) - 2\beta(n-1)}
\nonumber
\\[-8pt]
\\[-8pt]
\nonumber
&=&\frac{x}{\e^{2\beta(n-1)}L(\e^{2n-2})} \qquad\mbox{if } \e^{2\beta
(n-1)}L\bigl(\e^{2n-2}\bigr)
\leq x <\e^{2n\beta}L\bigl(\e^{2n}\bigr).
\end{eqnarray}
Thus, $\delta(x)$ has the required form (\ref{econditiondeltax})--(\ref{econditionax}) with
%
\begin{equation}
\label{edefinitionmainproofAkn} A_{k_{n}} = \e^{2\beta(n-1)}L\bigl(\e^{2n-2}\bigr)
\end{equation}
and
%
\begin{equation}
\label{edefinitionax} a(x) = \e^{2\beta(n-1)}L\bigl(\e^{2n-2}\bigr) =
A_{k_{n}}\qquad \mbox{if }A_{k_{n}}\leq x<A_{k_{n+1}}.
\end{equation}
Note also from (\ref{edeltaxfinalform}) that
\begin{eqnarray*}
1&\leq&\delta(x) <\frac{\e^{2\beta n}L(\e^{2n})}{\e^{2\beta
(n-1)}L(\e^{2n-2})} \\
&=& \e^{2\beta}\frac
{L(\e^{2n})}{L(\e^{-2}\e^{2n})}
\rightarrow \e^{2\beta} = c^{{1}/{\alpha
}},
\end{eqnarray*}
so that $\delta(x)\in[1, c^{{1}/{\alpha}} +\varepsilon]$ for large
enough $x$ when $\varepsilon> 0$ is fixed.

To complete step $1$, we need to prove that
$h_{R}(A_{k_{n}}x_{0})\rightarrow0$ as $n\rightarrow\infty$ for every
continuity point $x_{0}$ of $M_{R}(x)$. The discontinuity points of
$M_{R}$ are
%
\begin{equation}
\label{ediscontinuityMR} x = \e^{2k\beta},\qquad k\in\mathbb{Z}.
\end{equation}
To show $h_{R}(A_{k_{n}}x_{0})\rightarrow0$, note that, by Lemma~\ref
{lg1-tilde-g1}, it is enough to prove that $\tilde
{h}_{R}(A_{k_{n}}x_{0})\neq0$ for finitely many values of $n$, where
\[
\tilde{h}_{R}(x) = \e^{-\nu\lceil g_{1}(({1}/{(2\beta)})\operatorname{log
}x)\rceil_{\inplus}} - \e^{-\nu\lceil\tilde{g}_{1}(({1}/{(2\beta)})\operatorname
{log }x)\rceil_{\inplus}}.
\]
This holds only if for some integer $m\geq2$,
%
\begin{equation}
\label{eh{R}part} m +\log L\bigl(\e^{2m - 2}\bigr)\leq\frac{1}{2\beta}\log
A_{k_{n}}x_{0} < m +\log L\bigl(\e^{2m}\bigr).
\end{equation}
By Lemma~\ref{lfinite-number}, (\ref{eh{R}part}) holds for
infinitely many values of $n$ only if $x_{0} =\e^{2r\beta}$, $r\in\mathbb
{Z}$, which is a discontinuity point of $M_{R}(x)$ in (\ref
{ediscontinuityMR}).
Hence, $h_{R}(A_{k_{n}}x_{0})\rightarrow0$ as $n\rightarrow\infty$ for
every continuity point $x_{0}$ of $M_{R}(x)$.

\textit{Step} 2 (showing (\ref{ecorollary3,megyesi,lefttail})): In view of (\ref{erandomvariabledefinition}), we
are now interested in
%
\begin{equation}
\label{eFxlefttail} F_{-}(-x) = P \bigl(L \bigl(\e^{W_{q}}
\bigr)\e^{\beta W_{q}}(-1)^{W_{q}} < - x \bigr).
\end{equation}
Let $Z_{2} = \frac{W_{q}}{2}$ as in step 1. Note that (\ref
{eFxlefttail}) can be written as
%
\begin{eqnarray}
\label{elefttailexpressed} F_{-}(-x)&=&P \biggl(L \bigl(\e^{2Z_{2}}
\bigr)\e^{2\beta Z_{2}} > x, Z_{2}-\frac{1}{2}\mbox{ is integer}
\biggr)
\nonumber
\\
&=&P \biggl(L \bigl(\ee^{2(Z_{2}-{1}/{2})} \bigr)\e^{\beta}\e^{2\beta
(Z_{2}-{1}/{2})} >
x, Z_{2}-\frac{1}{2}\mbox{ is integer} \biggr)
\nonumber
\\[-8pt]
\\[-8pt]
\nonumber
&=&P \bigl(L\bigl(\ee^{2Z_{1}}\bigr)\e^{\beta}\e^{2\beta Z_{1}} > x
\bigr)
\\
&=&P \biggl(Z_{1} +\frac{1}{2} + \frac{1}{2\beta}\log L
\bigl(\e^{2Z_{1} + 1}\bigr) >\frac{1}{2\beta}\operatorname{log } x \biggr),\nonumber
\end{eqnarray}
where, in view of (\ref{eoddWq}),
%
\begin{equation}
\label{edefinitionZ1} P(Z_{1} \geq x) = h_{2}\bigl(\lceil x \rceil
\bigr)\e^{-\nu\lceil x \rceil}.
\end{equation}
Writing (\ref{elefttailexpressed}) as
\[
F_{-}(-x) = P \biggl(Z_{1} + \frac{1}{2\beta}\log L
\bigl(\e^{2Z_{1} + 1}\bigr) >\frac
{1}{2\beta}\operatorname{log } x -
\frac{1}{2} \biggr),
\]
the right-hand side has the form (\ref{eFxrighttail}) where
$L(\e^{2Z_{2}})$ is replaced by $L(\ee^{2Z_{1}})$ and $\frac{1}{2\beta
}\log x$ is replaced by $\frac{1}{2\beta}\log x - \frac{1}{2}$. Thus,
as in (\ref{ecompactrighttaill})--(\ref{edefinitiong0}), one
can write
%
\begin{equation}
\label{efirstg0tilde} F_{-}(-x) = P \biggl(Z_{1} \geq\tilde{g}
\biggl(\frac{1}{2\beta}\log x - \frac
{1}{2}\biggr) \biggr),
\end{equation}
where
%
\begin{equation}
\tilde{g}(y) = n\qquad \mbox{if } n-1 +\frac{1}{2\beta}\log L\bigl(\ee^{2n-2}
\bigr) \leq y < n +\frac{1}{2\beta}\log L\bigl(\ee^{2n}\bigr).
\end{equation}
The expression (\ref{efirstg0tilde}) can also be written as
%
\begin{equation}
\label{esecondg0tilde} F_{-}(-x) = P \biggl(Z_{1} \geq
\tilde{g}_{0}\biggl(\frac{1}{2\beta}\log x\biggr) \biggr),
\end{equation}
where $\tilde{g}_{0}(y) = \tilde{g}(y - \frac{1}{2})$ or, for $n\geq2$,
%
\begin{equation}
\tilde{g}_{0}(y) = n\qquad \mbox{if } n-\frac{1}{2} +
\frac{1}{2\beta
}\log L\bigl(\e^{2n - 1}\bigr) \leq y < n+\frac{1}{2}
+\frac{1}{2\beta} \log L\bigl(\e^{2n + 1}\bigr).
\end{equation}

We want to work with the intervals $[n - 1 +\frac{1}{2\beta}\log
L(\e^{2n - 2}), n +\frac{1}{2\beta}\log L(\e^{2n}))$ appearing in step 1,
and use the results of that step. Note that, on the interval $[n - 1
+\frac{1}{2\beta}\log L(\e^{2n - 2}), n +\frac{1}{2\beta}\log
L(\e^{2n}))$, the function $\tilde{g}_{0}$ has the form
%
\begin{equation}
\label{edefinitiontildeg0*} \tilde{g}_{0}(y) = %
\cases{ n-1, &\quad $\mbox{if
} \displaystyle n-1 +\frac{1}{2\beta}\log L\bigl(\e^{2n - 2}\bigr)\leq y < n-
\frac{1}{2} +\frac{1}{2\beta}\log L\bigl(\e^{2n - 1}\bigr) $,
\vspace*{2pt}
\cr
n, &\quad $\mbox{if } \displaystyle n-\frac{1}{2} +\frac{1}{2\beta}\log L
\bigl(\e^{2n - 1}\bigr) \leq y < n+\frac{1}{2\beta} \log L
\bigl(\e^{2n}\bigr)$.} %
\end{equation}
Defining
%
\begin{equation}
\label{edefinitionI0} I_{0}(y) = %
\cases{ -1, &\quad $\mbox{if }\displaystyle n-1 +
\frac{1}{2\beta}\log L\bigl(\e^{2n - 2}\bigr)\leq y < n-\frac{1}{2}
+\frac{1}{2\beta}\log L\bigl(\e^{2n - 1}\bigr) $,\vspace*{2pt}
\cr
0, & \quad$
\mbox{if } \displaystyle n-\frac{1}{2} +\frac{1}{2\beta}\log L\bigl(\e^{2n - 1}
\bigr) \leq y < n+\frac{1}{2\beta} \log L\bigl(\e^{2n}\bigr)$,}
\end{equation}
and combining (\ref{edefinitiong0}), (\ref{edefinitiontildeg0*}) and
(\ref{edefinitionI0}), we have
%
\begin{equation}
\label{erelationI0g0} \tilde{g}_{0}(y) = g_{0}(y) +
I_{0}(y).
\end{equation}

Continuing with (\ref{esecondg0tilde}), note further that, by
using (\ref{edefinitionZ1}) and (\ref{erelationI0g0}),
%
\begin{eqnarray}
\label{eexpandleftatil} F_{-}(- x)&=&h_{2} \biggl(
\tilde{g}_{0}\biggl(\frac{1}{2\beta}\log x\biggr)
\biggr)\e^{-\nu\tilde{g}_{0}(({1}/{(2\beta)})\log x)}
\nonumber
\\[-8pt]
\\[-8pt]
\nonumber
&=&\e^{-\nu I_{0}(({1}/{(2\beta)})\operatorname{log } x)}h_{2} \biggl(g_{0}\biggl(
\frac
{1}{2\beta}\operatorname{log } x\biggr) + I_{0}\biggl(
\frac{1}{2\beta}\operatorname{log } x\biggr) \biggr)\e^{-\nu g_{0}(({1}/{(2\beta)})\operatorname{log } x)}.
\end{eqnarray}
We want to write $F_{-}(- x)$ as in (\ref{ecorollary3,megyesi,lefttail}) of Theorem~\ref{tmain-semistable} (where by Lemma~\ref{ltwo-slow-var}, we can take a slowly varying function $l_{2}^{*}$
which is asymptotically equivalent to $l_{1}^{*}$). We need the
notation for the intervals appearing in (\ref
{edefinitiontildeg0*})--(\ref{edefinitionI0}), namely, for
$n\geq1$,
\begin{eqnarray*}
D_{n} &=& \biggl[n-1+\frac{1}{2\beta}\operatorname{log }L
\bigl(\e^{2n-2}\bigr), n-\frac
{1}{2}+\frac{1}{2\beta}
\operatorname{log }L\bigl(\e^{2n-1}\bigr)\biggr),
\\
E_{n} &=&\biggl [n-\frac{1}{2}+\frac{1}{2\beta}\operatorname{log
}L\bigl(\e^{2n-1}\bigr), n+\frac
{1}{2\beta}\operatorname{log }L
\bigl(\e^{2n}\bigr)\biggr).
\end{eqnarray*}
We also need a similar notation without the slowly varying function
$L$, that is, for $n\geq1$,
\[
D_{n}^{\prime} = \bigl[n-1, n-\tfrac{1}{2}\bigr),\qquad
E_{n}^{\prime} = \bigl[n-\tfrac{1}{2}, n\bigr).
\]
Set also
%
\begin{eqnarray}
\label{eunionfunction} D &=& \bigcup_{n=1}^{\infty}
D_{n},\qquad E=\bigcup_{n=1}^{\infty}
E_{n},
\nonumber
\\[-8pt]
\\[-8pt]
\nonumber
 D^{\prime} &=& \bigcup_{n=1}^{\infty}
D_{n}^{\prime},\qquad E^{\prime} = \bigcup
_{n=1}^{\infty}E_{n}^{\prime}.
\end{eqnarray}

As in (\ref{eexpressedrighttaill}), we can now write (\ref
{eexpandleftatil}) as
\begin{eqnarray*}
F_{-}(-x) &=& x^{-\alpha}\frac{h_{2}(g_{0}(({1}/{(2\beta)})\operatorname{log }
x) + I_{0}(({1}/{(2\beta)})\operatorname{log } x))}{c_{1}h_{1}(g_{0}(
({1}/{(2\beta)})\operatorname{log } x))}l_{1}^{*}(x)c_{1}\e^{-\nu I_{0}(
({1}/{(2\beta)})\log x)}\\
&&{}\times\e^{-\nu(\lceil g_{1}(({1}/{(2\beta)})\operatorname{log
}x)\rceil_{\inplus} - g_{1}(({1}/{(2\beta)})\operatorname{log }x))},
\end{eqnarray*}
where $\alpha= \frac{\nu}{2\beta}$ and $l_{1}^{*}(x)$ is given in (\ref
{el1*x}). This can also be written as
\[
F_{-}(-x) = x^{-\alpha}l_{2}^{*}(x)
\bigl(M_{L}\bigl(-\delta(x)\bigr) + h_{L}(x)\bigr),
\]
where
%
\begin{eqnarray}
l_{2}^{*}(x) &=& \frac{h_{2}(g_{0}(({1}/{(2\beta)})\operatorname{log } x) +
I_{0}(({1}/{(2\beta)})\operatorname{log } x))}{c_{1}h_{1}(g_{0}(({1}/{(2\beta)})\operatorname{log } x))} l_{1}^{*}(x),
\\
\label{edefinitionMLdeltax} M_{L}\bigl(-\delta(x)\bigr) &=& c_{1}
\e^{-\nu([{1}/{2} + \tilde{g}_{1}(({1}/{(2\beta)})\operatorname{log }x)] -
\tilde{g}_{1}(({1}/{(2\beta)})\operatorname{log }x))},
\\
\label{eequationhL} h_{L}(x)& =& c_{1}\e^{-\nu I_{0}(({1}/{(2\beta)})\log x)}\e^{-\nu(\lceil
g_{1}(({1}/{(2\beta)})\operatorname{log }x)\rceil_{\inplus} - g_{1}(
({1}/{(2\beta)})\operatorname{log }x))}
\nonumber
\\[-8pt]
\\[-8pt]
\nonumber
&&{}- c_{1} \e^{-\nu([{1}/{2} + \tilde
{g}_{1}(({1}/{(2\beta)})\operatorname{log }x)] - \tilde{g}_{1}(({1}/{(2\beta)})\operatorname{log }x))}.
\end{eqnarray}

By using (\ref{eh2h1ratio})--(\ref{eh1ratio}), we have
\[
\frac{h_{2}(g_{0}(({1}/{(2\beta)})\operatorname{log } x) + I_{0}
(({1}/{(2\beta)})\operatorname{log } x))}{c_{1}h_{1}(g_{0}
(({1}/{(2\beta)})\operatorname{log }
x))}\rightarrow1 \qquad\mbox{as } x\rightarrow\infty.
\]
Hence, $\frac{l_{2}^{*}(x)}{l_{1}^{*}(x)}\rightarrow1$, as
$x\rightarrow\infty$, that is, $l_{2}^{*}(x)$ and $l_{1}^{*}(x)$ are
two asymptotically equivalent functions. By the definition of $I_{0}$
and using Lemma~\ref{ll1*}, $l_{2}^{*}(x)$ is right-continuous and
slowly varying.

The function $\delta(x)$ appearing in (\ref
{edefinitionMLdeltax}) is the same as in (\ref{edefinitiondeltax})--(\ref{edeltaxfinalform}) of step 1, while the
function $M_{L}(-x)$ is defined as
%
\begin{equation}
M_{L}(-x) = c_{1}\e^{-\nu([{1}/{2}+({1}/{(2\beta)})\log x] -( {1}/{(2\beta)})\log x)},\qquad x > 0.
\end{equation}
It is left-continuous when $x>0$, and also periodic with multiplicative
period $\e^{2\beta}= c^{{1}/{\alpha}}$. Thus, $M_{L}(x)$ for $x <0$
is left-continuous as required in Theorem~\ref{tmain-semistable}. The
discontinuity points of $M_{L}(-x)$ are
%
\begin{equation}
\label{ediscontinuityML} x = \e^{\beta(2k +1 )},\qquad k\in\mathbb{Z}.
\end{equation}

To conclude the proof of step 2, we need to show that
$h_{L}(A_{k_{n}}x_{0})\rightarrow0$ as $n\rightarrow\infty$ for every
continuity point $x_{0}$ of $M_{L}(-x)$, that is, $x_{0}$ different
from (\ref{ediscontinuityML}). For this, we rewrite $h_{L}(x)$ as
follows. Observe that
\[
\e^{-\nu I_{0}(y)} = \e^{\nu}1_{D}(y) + 1_{E}(y)
\]
and
\[
\bigl(\e^{\nu}1_{D^{\prime}}(y) + 1_{E^{\prime}}(y)
\bigr)\e^{-\nu(\lceil y \rceil_{\inplus} -
y)} = \e^{-\nu([{1}/{2} + y] - y)},
\]
where after taking the logs, using $\lceil y \rceil_{\inplus} =[ y ] +
1$ and simplification, the last identity is equivalent to
$[y]1_{D^{\prime}}(y) + ([y] + 1)1_{E^{\prime}}(y) = [\frac{1}{2} + y]$ and can
be seen easily by drawing a picture. By using these identities and (\ref
{eequationhL}), we can write
\begin{eqnarray*}
c_{1}^{-1}h_{L}(x) &=& \biggl(\e^{\nu}1_{D}
\biggl(\frac{1}{2\beta}\log x\biggr) + 1_{E}\biggl(
\frac{1}{2\beta}\log x\biggr)\biggr) \e^{-\nu(\lceil g_{1}(({1}/{(2\beta)
})\operatorname{log }x)\rceil_{\inplus} - g_{1}(({1}/{(2\beta)})\operatorname{log
}x))}
 \\
 &&{}-\e^{-\nu([{1}/{2} +
\tilde{g}_{1}(({1}/{(2\beta)})\operatorname{log }x)] - \tilde{g}_{1}(
({1}/{(2\beta)})\operatorname{log }x))}
\nonumber
\\
&=& h_{1, L}(x)\e^{-\nu(\lceil g_{1}(({1}/{(2\beta)})\operatorname{log }x)\rceil
_{\inplus} - g_{1}(({1}/{(2\beta)})\operatorname{log }x))} + h_{2,
L}(x),
\end{eqnarray*}
where
\begin{eqnarray*}
h_{1, L}(x)& =& \e^{\nu}1_{D}\biggl(
\frac{1}{2\beta}\operatorname{log }x\biggr) + 1_{E}\biggl(
\frac
{1}{2\beta}\operatorname{log }x\biggr) - \e^{\nu}1_{D^{\prime}}
\biggl({g}_{1}\biggl(\frac{1}{2\beta
}\operatorname{log }x\biggr)
\biggr)\\
&&{} - 1_{E^{\prime}}\biggl({g}_{1}\biggl(\frac{1}{2\beta}
\operatorname{log }x\biggr)\biggr)
,
\\
h_{2, L}(x) &=& \e^{-\nu([{1}/{2} + {g}_{1}(({1}/{(2\beta)})\operatorname{log
}x)] - {g}_{1}(({1}/{(2\beta)})\operatorname{log }x))}\\
&&{}- \e^{-\nu([{1}/{2} +
\tilde{g}_{1}(({1}/{(2\beta)})\operatorname{log }x)] - \tilde{g}_{1}(
({1}/{(2\beta)})\operatorname{log }x))}.
\end{eqnarray*}
It is therefore enough to show that $h_{1,
L}(A_{k_{n}}x_{0})\rightarrow0$ and $h_{2,
L}(A_{k_{n}}x_{0})\rightarrow0$, as $n\rightarrow\infty$.
From~(\ref{edefinitiong1}), (\ref{edefinitiong1tilde}) and (\ref
{eunionfunction}), $h_{1, L}(A_{k_{n}}x_{0})\neq0$ if, for some
integer $m\geq1$,
%
\begin{equation}
\label{eh{L}firstpart} m - \frac{1}{2} +\log L\bigl(\e^{2m - 1}\bigr)\leq
\frac{1}{2\beta}\log A_{k_{n}}x_{0} < m - \frac{1}{2} +
\log L\bigl(\e^{2m}\bigr).
\end{equation}
(To see this, partition $[m - 1 + \frac{1}{2\beta}\log L(\e^{2m-2}), m +
\frac{1}{2\beta}\log L(\e^{2m}))$ into four subintervals $[m - 1 + \frac
{1}{2\beta}\log L(\e^{2m-2}), m - 1 + \frac{1}{2\beta}\log L(\e^{2m}))$,
$[m - 1 + \frac{1}{2\beta}\log L(\e^{2m}), m - \frac{1}{2} + \frac
{1}{2\beta}\log L(\e^{2m-1}))$, $[m - \frac{1}{2} + \frac{1}{2\beta}\log
L(\e^{2m-1}), m - \frac{1}{2} + \frac{1}{2\beta}\log L(\e^{2m}))$, $[m -
\frac{1}{2} + \frac{1}{2\beta}\log L(\e^{2m}), m + \frac{1}{2\beta}\log
L(\e^{2m}))$ and check that the function is nonzero only on the third
subinterval as given in (\ref{eh{L}firstpart}).)
By Lemma~\ref{lfinite-number}, (\ref{eh{L}firstpart}) holds for
infinitely many values of $n$ only if $x_{0} =\e^{\beta(2r + 1)}$ which
is a discontinuity point of $M_{L}(-x)$ in (\ref{ediscontinuityML}).
To show $h_{2, L}(A_{k_{n}}x_{0})\rightarrow0$, note that, by Lemma~\ref{lg1-tilde-g1}, it is enough to prove that $\tilde{h}_{2,
L}(A_{k_{n}}x_{0})\neq0$ for finitely many values of $n$, where
\[
\tilde{h}_{2, L}(x) = \e^{-\nu[{1}/{2} + {g}_{1}(({1}/{(2\beta)})\operatorname{log }x)]} - \e^{-\nu[{1}/{2} + \tilde{g}_{1}(({1}/{(2\beta)})\operatorname{log }x)]}.
\]
By using (\ref{edefinitiong1}) and (\ref{edefinitiong1tilde}), the
relation $\tilde{h}_{2, L}(A_{k_{n}}x_{0}) = 0$ holds only if, for some
integer $m\geq1$,
%
\begin{equation}
\label{eh{L}secondpart} m - \frac{1}{2} +\log L\bigl(\e^{2m - 2}\bigr)\leq
\frac{1}{2\beta}\log A_{k_{n}}x_{0} < m - \frac{1}{2} +
\log L\bigl(\e^{2m}\bigr).
\end{equation}
(To see this, draw a plot of $g_{1}(y)$ and $\tilde{g}_{1}(y)$ for $y$
in $[m - 1+\frac{1}{2\beta}\log L(\e^{2m-2}), m - \frac{1}{2}\log
L(\e^{2m}))$, and note that $\tilde{g}_{1}(y) = m - \frac{1}{2}$ at $y =
m - \frac{1}{2}+\frac{1}{2\beta}\log L(\e^{2m-2})$ and $g_{1}(y) = m -
\frac{1}{2}+\frac{1}{2\beta}\log L(\e^{2m})$.)
By Lemma~\ref{lfinite-number}, (\ref{eh{L}secondpart}) holds for
infinitely many values of $n$ only if $x_{0} =\e^{\beta(2r + 1)}$ which
is a discontinuity point of $M_{L}(-x)$ in (\ref{ediscontinuityML}).
Hence, $h_{L}(A_{k_{n}}x_{0})\rightarrow0$ as $n\rightarrow\infty$ for
every continuity point $x_{0}$ of $M_{L}(-x)$.

\textit{Step} 3 (Deriving subsequence $k_{n}$): We conclude
the proof of the theorem by showing that $k_{n}$ is given by (\ref
{edefinitionknAknBkn}). In view of the discussion following Theorem~\ref{tmain-semistable}, we want to choose a sequence $\tilde{A}_{n}$
satisfying (\ref{edefineAntilde})--(\ref{econditionAntilde})
such that $k_{n}$ given by (\ref{edefinitionknAknBkn}) now satisfies
(\ref{edefinitionAk{n}}).
We define such sequence $\tilde{A}_{n}$ as
%
\begin{eqnarray}
\operatorname{log }\tilde{A}_{n} &=& 2\beta(m-1)+ \log L
\bigl(\e^{2m-2}\bigr)
\nonumber
\\[-8pt]
\\[-8pt]
\nonumber
&&{} + \frac
{{(\operatorname{log } n - \operatorname{log } k_{m})(2\beta+ \operatorname{log } L(\e^{2m}) -
\operatorname{log }L(\e^{2m-2}))}}{\operatorname{log }k_{m+1} - \operatorname{log
}k_{m}}
\qquad \mbox{if }k_{m} \leq n < k_{m+1}, m\geq1.
\end{eqnarray}
The sequence $\tilde{A}_{n}$ satisfies (\ref{econditionAntilde}). If
$k_{m} \leq n < k_{m+1} - 1$,
the last limit in (\ref{econditionAntilde}) follows from
\[
\operatorname{log }\tilde{A}_{n+1} - \operatorname{log }
\tilde{A}_{n} = \frac{{(\operatorname
{log } n - \operatorname{log } (n +1))(2\beta+ \operatorname{log } L(\e^{2m}) - \operatorname
{log }L(\e^{2m-2}))}}{\operatorname{log }k_{m+1} - \operatorname{log }k_{m}}\rightarrow 0.
\]
If $n = k_{m+1} - 1$, the limit follows from
\begin{eqnarray*}
&&\operatorname{log }\tilde{A}_{n+1} - \operatorname{log }
\tilde{A}_{n}\nonumber\\
&&\quad= 2\beta+ \log L\bigl(\e^{2m}\bigr) - \log L
\bigl(\e^{2m-2}\bigr)
\\
&&\qquad{} -\frac{{(\operatorname{log }(k_{m+1} - 1) - \operatorname{log }
k_{m})(2\beta+ \operatorname{log } L(\e^{2m}) - \operatorname{log }L(\e^{2m-2}))}}{\operatorname
{log }k_{m+1} - \operatorname{log }k_{m}}\rightarrow0
\nonumber
\end{eqnarray*}
since $\log L(\e^{2m}) - \log L(\e^{2m-2})\rightarrow0$, and
\[
\frac{{\operatorname{log }(k_{m+1} - 1) - \operatorname{log } k_{m}}}{\operatorname{log
}k_{m+1} - \operatorname{log }k_{m}}\rightarrow1.
\]

Next we show (\ref{edefineAntilde}), that is, $n\tilde
{A}_{n}^{-\alpha}l_{1}^{*}(\tilde{A}_{n})\rightarrow1$, as
$n\rightarrow\infty$, where $\alpha=\frac{\nu}{2\beta}$ and
$l_{1}^{*}$ is as defined in (\ref{el1*x}). When $k_{m} \leq n <
k_{m+1}$, observe that
%
\begin{eqnarray}
\label{easymptoticsAntilde} \log n\tilde{A}_{n}^{-\alpha}l_{1}^{*}(
\tilde{A}_{n}) &=&\log\frac{n
l_{1}^{*}(\tilde{A}_{n})}{\tilde{A}_{n}^{\nu/2\beta}}
\nonumber
\\
&=&\log\frac{n l_{1}^{*}(\tilde{A}_{n})}{\e^{(m-1)\nu}L(\e^{2m-2})^{\nu
/2\beta}}\nonumber\\
&&{}+ \frac{\nu+ ({\nu}/{(2\beta)})\operatorname{log } L(\e^{2m}) -
({\nu}/{(2\beta)})\operatorname{log }L(\e^{2m-2})}{\operatorname{log }k_{m+1} - \operatorname{log
}k_{m}}\log \biggl(\frac{k_{m}}{n}
\biggr)\qquad
\\
&\sim&\log n + \log\frac{l_{1}^{*}(\tilde
{A}_{n})}{h_{1}(m-1)L(\e^{2m-2})^{\nu/2\beta}}- \log{k_{m}}
\nonumber
\\
&&{} +\frac{\nu+ ({\nu}/{(2\beta)})\operatorname{log } L(\e^{2m}) -
({\nu}/{(2\beta)})\operatorname{log }L(\e^{2m-2})}{\operatorname{log }k_{m+1} - \operatorname{log
}k_{m}}\log \biggl(\frac{k_{m}}{n} \biggr).\nonumber
\end{eqnarray}
Now observe that as $n\rightarrow\infty$, we have $m\rightarrow\infty$,
and thus $\frac{k_{m}}{n}$ is bounded and
\[
\frac{\nu+ ({\nu}/{(2\beta)})\operatorname{log } L(\e^{2m}) -( {\nu}/{(2\beta)
})\operatorname{log }L(\e^{2m-2})}{\operatorname{log }k_{m+1} - \operatorname{log
}k_{m}}\rightarrow1.
\]
Thus, (\ref{easymptoticsAntilde}) is asymptotically equivalent to
%
\begin{equation}
\label{epseudosimplified1} \log\frac{l_{1}^{*}(\tilde{A}_{n})}{h_{1}(m-1)L(\e^{2m-2})^{\nu/2\beta}}.
\end{equation}
By the relation (\ref{easymptoticequivalent}) in the \hyperref[sauxiliary]{Appendix}, $l_{1}^{*}(\tilde{A}_{n})\sim h_{1}(g_{2}(\frac{1}{2\beta
}\operatorname{log }\tilde{A}_{n}))\e^{\nu\tilde{g}_{1}^{*}(({1}/{(2\beta)})\operatorname
{log }\tilde{A}_{n})}$ and hence (\ref{easymptoticsAntilde}) is also
asymptotically equivalent to
%
\begin{equation}
\label{epseudosimplified} \log\frac{h_{1}(g_{2}(({1}/{(2\beta)})\operatorname{log }\tilde{A}_{n}))\e^{\nu
\tilde{g}_{1}^{*}(({1}/{(2\beta)})\operatorname{log }\tilde
{A}_{n})}}{h_{1}(m-1)L(\e^{2m-2})^{\nu/2\beta}}.
\end{equation}
Since $k_{m} \leq n < k_{m+1}$, we have
\[
2\beta(m-1) + \log L\bigl(\e^{2m-2}\bigr) \leq\log\tilde{A}_{n}
< 2\beta m + \log L\bigl(\e^{2m}\bigr)
\]
and, by (\ref{etildeg1*}), $\frac{\e^{\nu\tilde{g}_{1}^{*}(
({1}/{(2\beta)})\operatorname{log}\tilde{A}_{n})}}{L(\e^{2m-2})^{\nu/2\beta}} = 1$.
Hence, (\ref{epseudosimplified}) simplifies to $\log\frac
{h_{1}(m-1+\kappa)}{h_{1}(m-1)},\mbox{where }0\leq\kappa< 1$. But as
$n\rightarrow\infty$, we have $m\rightarrow\infty$ and thus $\frac
{h_{1}(m-1+\kappa)}{h_{1}(m-1)}\rightarrow1$ by using (\ref{eh1ratio}).
This proves that $\log n\tilde{A}_{n}^{-\alpha}l_{1}^{*}(\tilde
{A}_{n})\rightarrow0$ and thus $n\tilde{A}_{n}^{-\alpha
}l_{1}^{*}(\tilde{A}_{n})\rightarrow1$, as $n\rightarrow\infty$.

Finally, we show that $k_{n}$ defined in (\ref
{edefinitionknAknBkn}) satisfies (\ref{edefinitionAk{n}}). Define
$a_{n} = A_{k_{n}} = \e^{2\beta(n-1)}L(\e^{2n-2})$.
Hence,
\[
\log a_{n} = \log A_{k_{n}} = 2\beta(n-1) + \log L
\bigl(\e^{2n-2}\bigr).
\]
Now observe that $\tilde{A}_{k_{n}} = a_{n}$ and thus (\ref
{edefinitionAk{n}}) is satisfied.
\end{pf}

The partial sums (\ref{epartialsum}) involve centering constants
$B_{k_{n}}$ defined in (\ref{edefiBkn}). As in the stable case, one
can expect to replace $B_{k_{n}}$ by $k_{n}E X$ when $1<\alpha<2$, and
to show the convergence of (\ref{epartialsum}) without $B_{k_{n}}$
when $0<\alpha<1$. The next result shows that this is indeed the case.

\begin{proposition} \label{pcentering} Suppose that the assumptions of
Theorem~\ref{tmain-th} hold. Let
%
\begin{eqnarray}
\label{eeta1} \zeta&=& - \frac{1 - \e^{-\nu}}{1 - \e^{2\beta- \nu}} -
\e^{\beta(2\lceil
({1}/{\nu})\log c_{1}\rceil- 1)}
\bigl(c_{1}\e^{-\nu(\lceil({1}/{\nu
})\log c_{1}\rceil- 1)} - 1\bigr)
\nonumber
\\[-8pt]
\\[-8pt]
\nonumber
& &{}+ c_{1}\frac{(1- \e^{-\nu})\e^{\nu- \beta}}{1 - \e^{2\beta- \nu
}}\e^{(2\beta- \nu)\lceil({1}/{\nu})\log c_{1}\rceil}.
\end{eqnarray}
If $0<\alpha<1$, then
\[
\frac{B_{k_{n}}}{A_{k_{n}}}\rightarrow\zeta, \qquad\frac
{1}{A_{k_{n}}}\sum
_{j=1}^{k_{n}} X_{j}\stackrel{d}\rightarrow
Y+\zeta
\]
and if $1<\alpha<2$, then
\[
\frac{k_{n}E X - B_{k_{n}}}{A_{k_{n}}}\rightarrow- \zeta,\qquad \frac
{1}{A_{k_{n}}} \Biggl\{\sum
_{j=1}^{k_{n}}X_{j} - k_{n}E {X}
\Biggr\} \stackrel{d}\rightarrow Y + \zeta,
\]
where $Y$ follows the semi-stable law characterized by (\ref
{echarcterization0}) and (\ref{echarcterization}).
\end{proposition}

\begin{pf}
\textit{Case} $0<\alpha<1$: It is enough to show the
convergence of $\frac{B_{k_{n}}}{A_{k_{n}}} = \frac
{k_{n}}{A_{k_{n}}}\int_{{1}/{k_{n}}}^{1 - {1}/{k_{n}}} Q(s) \,\mathrm{d}s$
to~$\zeta$, where $Q(s)$ is defined in (\ref{edefinitionQs}). For
fixed $s_{1}$ and $s_{2}$, write
%
\begin{eqnarray}
\label{eintegralrepre} &&\frac{k_{n}}{A_{k_{n}}}\int_{{1}/{k_{n}}}^{1 - {1}/{k_{n}}}
Q(s) \,\mathrm{d}s
\nonumber
\\[-8pt]
\\[-8pt]
\nonumber
&&\quad = \frac{k_{n}}{A_{k_{n}}}\int_{{1}/{k_{n}}}^{s_{1}} Q(s)
\,\mathrm{d}s + \frac{k_{n}}{A_{k_{n}}}\int_{s_{1}}^{s_{2}} Q(s) \,\mathrm{d}s +
\frac
{k_{n}}{A_{k_{n}}}\int_{s_{2}}^{1 - {1}/{k_{n}}} Q(s) \,\mathrm{d}s.
\end{eqnarray}
Observe first that, for fixed $s_{1}$ and $s_{2}$, the second term in
(\ref{eintegralrepre}) converges to zero. Indeed, this follows from
the fact that $\frac{k_{n}}{A_{k_{n}}}\rightarrow0$. For the latter
convergence, note from (\ref{edefinitionknAknBkn}) that
%
\begin{equation}
\label{einterepresecond} \frac{k_{n}}{A_{k_{n}}}\sim\frac{\e^{(n-1)\nu}}{h_{1}(n-1)}\frac
{1}{L(\e^{2n-2})\e^{2\beta(n-1)}}.
\end{equation}
For arbitrarily small $\delta>0$, by using Potter's bounds for $L$ and
Lemma~\ref{lh1-bounds} for $h_{1}$, the right-hand side of (\ref
{einterepresecond}) is bounded by $C\mathrm{e}^{(\nu- 2\beta+ \delta
)(n-1)}\rightarrow0$, as long as $\nu- 2\beta+ \delta< 0$.

Consider now the third term in (\ref{eintegralrepre}), involving
the function $Q(s)$ for values of $s$ close to $1$. The function $Q(s)$
is defined as the inverse of the distribution function $F(x) =
P(L(\e^{W_{q}})\e^{\beta W_{q}}(-1)^{W_{q}}\leq x)$. Since we are
interested in $Q(s)$ for $s$ close to $1$, it is enough to look at the
function for $x>0$. For $x>0$, the function $F(x)$ has jumps at points
$x = L(\e^{2n})\e^{2\beta n}$ of size
\[
P(W_q = 2n) = P\biggl(\frac{W_{q}}{2}\geq n, W_{q}
\mbox{ is even}\biggr) - P\biggl(\frac{W_{q}}{2}\geq n+1, W_{q}
\mbox{ is even}\biggr).
\]
This means that, for $s$ close to $1$, the inverse function $Q(s)$ has
jumps at points $s = 1 - P(\frac{W_{q}}{2}\geq n, W_{q}\mbox{ is
even})$ of size $L(\e^{2n})\e^{2\beta n} - L(\e^{2n-2})\e^{2\beta(n-1)}$.
Moreover, $Q(s) = L(\e^{2n})\e^{2\beta n}$ when $1 - P(\frac
{W_{q}}{2}\geq n, W_{q}\mbox{ is even}) \leq s < 1 - P(\frac
{W_{q}}{2}\geq n + 1, W_{q}\mbox{ is even})$. (If this step is unclear,
the reader may want to draw a picture.) Note that the jump points satisfy
\[
1 - s = P\biggl(\frac{W_{q}}{2}\geq n, W_{q}\mbox{ is even}
\biggr) = h_{1}(n)\e^{-\nu n}
\]
by (\ref{eevenWq}).

Assuming for simplicity that $\frac{\e^{\nu(n-1)}}{h_{1}(n-1)}$
are integers so that $k_{n} = \frac{\e^{\nu(n-1)}}{h_{1}(n-1)}$ and
taking $s_{2} = 1 - h_{1}(n_{1})\e^{-\nu n_{1}}$, we can write,
\begin{eqnarray*}
\label{erighttailcenterapproximation} &&\frac{k_{n}}{A_{k_{n}}}\int_{s_{2}}^{1 - h_{1}(n-1)\e^{-\nu(n-1)}}
Q(s) \,\mathrm{d}s
\\
&&\quad=\frac{k_{n}}{A_{k_{n}}}\sum_{m=n_{1}}^{n-2}L
\bigl(\e^{2m}\bigr)\e^{2\beta
m}\bigl(h_{1}(m)\e^{-\nu m}
- h_{1}(m+1)\e^{-\nu(m+1)}\bigr)
\\
&&\quad= \frac{\e^{\nu(n-1)}}{h_{1}(n-1)\e^{2\beta(n-1)}L(\e^{2n-2})}\sum_{m=n_{1}}^{n-2}L
\bigl(\e^{2m}\bigr)\e^{2\beta m}h_{1}(m)\e^{-\nu m}
\biggl(1 - \frac
{h_{1}(m+1)}{h_{1}(m)}\e^{-\nu}\biggr) \\
&&\quad=: I_{1} +
I_{2},
\end{eqnarray*}
where, for fixed $K$,
\begin{eqnarray*}
I_{1}&=&\frac{\e^{\nu(n-1)}}{h_{1}(n-1)\e^{2\beta(n-1)}L(\e^{2n-2})}\sum_{m=n_{1}}^{n-K}L
\bigl(\e^{2m}\bigr)\e^{2\beta m}h_{1}(m)\e^{-\nu m}
\biggl(1 - \frac
{h_{1}(m+1)}{h_{1}(m)}\e^{-\nu}\biggr),
\\
I_{2}& =& \frac{\e^{\nu(n-1)}}{h_{1}(n-1)\e^{2\beta(n-1)}L(\e^{2n-2})}\sum_{m=n-K}^{n-2}L
\bigl(\e^{2m}\bigr)\e^{2\beta m}h_{1}(m)\e^{-\nu m}
\biggl(1 - \frac
{h_{1}(m+1)}{h_{1}(m)}\e^{-\nu}\biggr).
\end{eqnarray*}
For the term $I_{2}$, note that, after changing $m$ to $n - j$ in the sum,
\[
I_{2} = \e^{2\beta- \nu}\sum_{j=2}^{K}
\frac
{L(\e^{2(n-j)})}{L(\e^{2(n-1)})}\e^{-(2\beta- \nu)j}\frac
{h_{1}(n-j)}{h_{1}(n-1)}\biggl(1 -
\frac{h_{1}(n - j+1)}{h_{1}(n - j)}\e^{-\nu
}\biggr).
\]
By using (\ref{eh1ratio}), we get that
%
\begin{equation}
\label{efirsttermcentering} I_{2}\rightarrow \e^{2\beta-\nu}
\bigl(1-\e^{-\nu}\bigr)\sum_{j=2}^{K}\e^{-(2\beta-
\nu)j}
= \bigl(1 - \e^{-\nu}\bigr)\frac{\e^{\nu- 2\beta}}{1 - \e^{\nu- 2\beta}}\bigl(1 - \e^{-(K-1)(2\beta- \nu)}
\bigr),
\end{equation}
as $n\rightarrow\infty$.
For the term $I_{1}$, we have similarly
\[
I_{1} = \e^{2\beta- \nu}\sum_{j=K}^{n-n_{1}}
\frac
{L(\e^{2(n-j)})}{L(\e^{2(n-2)})}\e^{-(2\beta- \nu)j}\frac
{h_{1}(n-j)}{h_{1}(n-1)}\biggl(1 -
\frac{h_{1}(n - j+1)}{h_{1}(n - j)}\e^{-\nu
}\biggr).
\]
For arbitrarily small $\delta>0$, by using Potter's bounds and Lemma~\ref{lh1-bounds}, we can write
%
\begin{equation}
\label{esecondtermcentering} |I_{1}| \leq C\sum_{j= K}^{n - n_{1}}
\e^{-(2\beta- \nu-\delta)j}.
\end{equation}
When $2\beta- \nu-\delta> 0$, the last bound is arbitrarily small
for large enough $K$. Together with~(\ref{efirsttermcentering}), this
shows that
\[
\frac{k_{n}}{A_{k_{n}}}\int_{s_{2}}^{1 - h_{1}(n-1)\e^{-(n-1)\nu}} Q(s) \,\mathrm{d}s =
I_{1} + I_{2}\rightarrow\bigl(1 - \e^{-\nu}\bigr)
\frac{\e^{\nu- 2\beta}}{1 -
\e^{\nu- 2\beta}} = - \frac{1 - \e^{-\nu}}{1 - \e^{2\beta- \nu
}},
\]
as $n\rightarrow\infty$.

Consider now the first term in (\ref{eintegralrepre}), involving
the function $Q(s)$ for values of $s$ close to 0. Here we need to
examine the function $F(x) = P(L(\e^{W_{q}})\e^{\beta
W_{q}}(-1)^{W_{q}}\leq x)$ for $x<0$. For $x<0$, the function $F(x)$
has jumps at $x = - L(\e^{2n+1})\e^{\beta(2n+1)}$ of size
\[
P(W_q = 2n+1) = P\biggl(\frac{W_{q} - 1}{2}\geq n,
W_{q}\mbox{ is odd}\biggr) - P\biggl(\frac{W_{q} -
1}{2}\geq n+1,
W_{q}\mbox{ is odd}\biggr).
\]
Moreover, $Q(s) =-L(\e^{2n+1})\e^{\beta(2n+ 1)}$ when $ P(\frac{W_{q} -
1}{2}\geq n+1, W_{q}\mbox{ is odd}) < s \leq P(\frac{W_{q} - 1}{2}\geq
n, W_{q}\mbox{ is odd})$. Note that, by (\ref{eoddWq}), the jump
points satisfy
\[
s = P\biggl(\frac{W_{q} - 1}{2}\geq n, W_{q}\mbox{ is odd}\biggr) =
h_{2}(n)\e^{-\nu
n}.
\]

Write the first term in (\ref{eintegralrepre}) as
%
\begin{equation}
\label{emixturetails} \frac{k_{n}}{A_{k_{n}}}\int_{h_{1}(n-1)\e^{-\nu
(n-1)}}^{h_{2}(l(n)-1)\e^{-\nu(l(n)-1)}}
Q(s) \,\mathrm{d}s +\frac
{k_{n}}{A_{k_{n}}}\int_{h_{2}(l(n)-1)\e^{-\nu(l(n)-1)}}^{s_{1}} Q(s) \,\mathrm{d}s
=: I_{1}^{*} + I_{2}^{*},
\end{equation}
where $l(n)$ is the integer such that
\[
\label{elefttailapprox} h_{2}\bigl(l(n)\bigr)\e^{-\nu l(n)}\leq
h_{1}(n-1)\e^{-\nu(n-1)} < h_{2}\bigl(l(n) - 1
\bigr)\e^{-\nu(l(n) -1)}
\]
or
\[
\label{einequalityofh1ln} h_{2}\bigl(l(n)\bigr)\e^{-\nu l(n)}\leq
h_{2}(n-1)\e^{-\nu(n-1 +({1}/{\nu})\log
({h_{2}(n-1)}/{(h_{1}(n-1))}))} < h_{2}\bigl(l(n) - 1
\bigr)\e^{-\nu(l(n)
-1)}.
\]
Note that, when $\frac{h_{2}(x)}{h_{1}(x)}\rightarrow c_{1}$ and $\frac
{1}{\nu}\log c_{1}$ is not an integer, or when $\frac{1}{\nu}\log
c_{1}$ is an integer and $\frac{h_{2}(x)}{h_{1}(x)}\uparrow c_{1}$, for
large values of $n$ one can take $l(n) = n-1 + \lceil\frac{1}{\nu}\log
c_{1}\rceil$. Indeed, this follows from
%
\begin{equation}
\label{econditionapprox1} \e^{-\nu}< \e^{-\nu(\lceil({1}/{\nu})\log
({h_{2}(n-1)}/{(h_{1}(n-1))})\rceil- ({1}/{\nu})\log
({h_{2}(n-1)}/{(h_{1}(n-1))}))}\leq1
\end{equation}
and the fact that
%
\begin{equation}
\label{econditionapprox2} \frac{h_{2}(n - 1 +\lceil({1}/{\nu})\log
({h_{2}(n-1)}/{(h_{1}(n-1))})\rceil)}{h_{2}(n-1)}\rightarrow1,
\end{equation}
as $n\rightarrow\infty$.

Now, taking $s_{1} = h_{2}(n_{2})\e^{-\nu n_{2}}$, we can write
$I_{2}^{*}$ in (\ref{emixturetails}) as
\[
I^{*}_{2} = - \frac{k_{n}}{A_{k_{n}}} \sum
_{m= n_{2}}^{l(n)-2} L\bigl(\e^{2m+1}
\bigr)\e^{\beta(2m+1)}\bigl(h_{2}(m)\e^{-\nu m} -
h_{2}(m+1)\e^{-\nu
(m+1)}\bigr).
\]

Following a similar calculation as done for the third term in (\ref
{eintegralrepre}), we get, as $n\rightarrow\infty$,
\begin{eqnarray*}
I^{*}_{2}&\rightarrow&- c_{1}
\frac{(1- \e^{-\nu})\e^{2(\nu- 2\beta
)}\e^{\beta}}{1 - \e^{\nu-2\beta}}\e^{-(\nu- 2\beta)\lceil({1}/{\nu
})\log c_{1}\rceil}\\
& =& c_{1}\frac{(1- \e^{-\nu})\e^{\nu- \beta}}{1 -
\e^{2\beta- \nu}}\e^{(2\beta- \nu)\lceil({1}/{\nu})\log c_{1}\rceil
}.
\end{eqnarray*}
One can write $I_{1}^{*}$ in (\ref{emixturetails}) as
\begin{eqnarray*}
I_{1}^{*} &=& -\frac{\e^{\nu(n-1)}}{h_{1}(n-1)\e^{2\beta
(n-1)}L(\e^{2n-2})}L\bigl(\e^{2l(n) -1}
\bigr)\e^{\beta(2l(n) - 1)}
\\
&&\hspace*{6pt}{}\times \bigl(h_{2}\bigl(l(n) - 1\bigr)\e^{-\nu(l(n)-1)} -
h_{1}(n-1)\e^{-\nu
(n-1)}\bigr)
\\
&=&-\frac{L(\e^{2l(n)-1})\e^{\beta(2l(n)-1)}}{L(\e^{2n-2})\e^{2\beta
(n-1)}}\biggl(\frac{h_{2}(l(n)-1)}{h_{1}(n-1)}\e^{-\nu(l(n)-n)}-1\biggr)
\\
&=&-\frac{L(\e^{2l(n)-1})}{L(\e^{2n-2})}\e^{\beta(2\lceil({1}/{\nu})\log
c_{1}\rceil-1)}\biggl(\frac{h_{2}(l(n)-1)}{h_{1}(n-1)}\e^{-\nu(\lceil
({1}/{\nu})\log c_{1}\rceil-1)}-1
\biggr)
.
\end{eqnarray*}
Now, by using (\ref{eh2h1ratio}) and (\ref{eh1ratio}), it can be seen that
\[
I_{1}^{*}\rightarrow-\e^{\beta(2\lceil({1}/{\nu})\log c_{1}\rceil-
1)}
\bigl(c_{1}\e^{-\nu(\lceil({1}/{\nu})\log c_{1}\rceil- 1)} - 1\bigr),
\]
as $n\rightarrow\infty$.

Now we consider the case when $\frac{1}{\nu}\log c_{1}$ is an integer
and $\frac{h_{2}(x)}{h_{1}(x)}\downarrow c_{1}$. We want to find $l(n)$
such that (\ref{einequalityofh1ln}) holds. Hence, we want
\[
\lim_{n\rightarrow\infty}\frac{h_{2}(n-1)}{h_{2}(l(n))} \e^{-\nu(n-1 +
({1}/{\nu})\log({h_{2}(n-1)}/{(h_{1}(n-1))}) - l(n))}\geq1.
\]
Take $l(n) = n - 2 +\lceil\frac{1}{\nu}\log\frac
{h_{2}(n-1)}{h_{1}(n-1)}\rceil$. Then, $\lim_{n\rightarrow\infty}\frac
{h_{2}(n-1)}{h_{1}(l(n))}\rightarrow1$. Now,
\begin{eqnarray*}
&&\lim_{n\rightarrow\infty}\e^{-\nu(n-1 + ({1}/{\nu})\log
({h_{2}(n-1)}/{(h_{1}(n-1))} )- n + 2 - \lceil({1}/{\nu})\log
({h_{2}(n-1)}/{(h_{1}(n-1)}))\rceil)} \\
&&\quad= \e^{-\nu(1+ ({1}/{\nu})\log c_{1} -
({1}/{\nu})\log c_{1} - 1)} =
\e^{0} = 1.
\end{eqnarray*}
We also need
\[
\frac{h_{2}(l(n)-1)}{h_{2}(n-1)} \e^{-\nu(l(n) - 1 - n + 1 - ({1}/{\nu})
\log({h_{2}(n-1)}/{(h_{1}(n-1))}))}> 1
\]
for large $n$.
For this, observe that $\frac{h_{2}(l(n)-1)}{h_{2}(n-1)}\rightarrow1$ and
\begin{eqnarray*}
&&\lim_{n\rightarrow\infty}\e^{-\nu(l(n) - 1 - n + 1 - ({1}/{\nu})\log
({h_{2}(n-1)}/{(h_{1}(n-1))}))} \\
&&\quad= \lim_{n\rightarrow\infty}\e^{-\nu(n- 3
+ \lceil({1}/{\nu})\log({h_{2}(n-1)}/{(h_{1}(n-1))})\rceil- n + 1 -
({1}/{\nu})\log({h_{2}(n-1)}/{(h_{1}(n-1))}))}
\\
&&\quad =\lim_{n\rightarrow\infty}\e^{-\nu(-2 + 1 + ({1}/{\nu})\log c_{1} -
({1}/{\nu})\log c_{1})} = \e^{-\nu}.
\end{eqnarray*}
Hence, when $\frac{1}{\nu}\log c_{1}$ is an integer, we have $\frac
{h_{2}(x)}{h_{1}(x)}\downarrow c_{1}$ and $\lceil\frac{1}{\nu}\log\frac
{h_{2}(x)}{h_{1}(x)}\rceil\downarrow\frac{1}{\nu}\log c_{1} + 1$, and
as in the previous calculations,
\[
I_{1}^{*}\rightarrow-\e^{\beta(2\lceil({1}/{\nu})\log c_{1}\rceil-
1)}
\bigl(c_{1}\e^{-\nu(\lceil({1}/{\nu})\log c_{1}\rceil- 1)} - 1\bigr)
\]
and
\[
I^{*}_{2}\rightarrow c_{1}\frac{(1- \e^{-\nu})\e^{\nu- \beta}}{1 -
\e^{2\beta- \nu}}\e^{(2\beta- \nu)\lceil({1}/{\nu})\log c_{1}\rceil
}.
\]

Finally, gathering the results above, we deduce the convergence to
the constant $\zeta$ given by~(\ref{eeta1}).

\textit{Case} $1<\alpha<2$: It is enough to show the
convergence of $\frac{k_{n}E X - B_{k_{n}}}{A_{k_{n}}}$ to $-\zeta$.
Using the fact that $E X = \int_{0}^{1} Q(s) \,\mathrm{d}s$, observe that
\[
\frac{k_{n}E X - B_{k_{n}}}{A_{k_{n}}} =\frac{k_{n}}{A_{k_{n}}}\int_{0}^{{1}/{k_{n}}}
Q(s) \,\mathrm{d}s +\frac{k_{n}}{A_{k_{n}}}\int_{1 -
{1}/{k_{n}}}^{1} Q(s) \,\mathrm{d}s.
\]
For simplicity, we assume that $\frac{\e^{(n-1)\nu}}{h_{1}(n-1)}$ is an
integer. To evaluate $\int_{1 - {1}/{k_{n}}}^{1} Q(s) \,\mathrm{d}s$, one
follows a similar procedure as in the case $0<\alpha<1$ to obtain
\begin{eqnarray*}
\label{eapproxcase2} &&\frac{k_{n}}{A_{k_{n}}}\int_{1 - h_{1}(n-1)\e^{-\nu(n-1)}}^{1}
Q(s) \,\mathrm{d}s
\\
&&\quad=\frac{k_{n}}{A_{k_{n}}}\sum_{m=n-1}^{\infty} L
\bigl(\e^{2m}\bigr)\e^{2\beta
m}\bigl(h_{1}(m)\e^{-\nu m}
- h_{1}(m+1)\e^{-\nu(m+1)}\bigr) =: \tilde{I}_{1} +
\tilde{I}_{2},
\end{eqnarray*}
where, for fixed $K$,
\begin{eqnarray*}
\tilde{I}_{1}&=&\frac{\e^{\nu(n-1)}}{h_{1}(n-1)\e^{2\beta
(n-1)}L(\e^{2n-2})}\sum_{m=n-1}^{n + K}L
\bigl(\e^{2m}\bigr)\e^{2\beta
m}h_{1}(m)\e^{-\nu m}
\biggl(1 - \frac{h_{1}(m+1)}{h_{1}(m)}\e^{-\nu}\biggr),
\\
\tilde{I}_{2} &= &\frac{\e^{\nu(n-1)}}{h_{1}(n-1)\e^{2\beta
(n-1)}L(\e^{2n-2})}\sum_{m=n + K}^{\infty}L
\bigl(\e^{2m}\bigr)\e^{2\beta
m}h_{1}(m)\e^{-\nu m}
\biggl(1 - \frac{h_{1}(m+1)}{h_{1}(m)}\e^{-\nu}\biggr).
\end{eqnarray*}
Similar to the case $0<\alpha<1$, one can show that
\[
\frac{k_{n}}{A_{k_{n}}}\int_{1 - h_{1}(n-1)\e^{-(n-1)\nu}}^{1} Q(s) \,\mathrm{d}s =
\tilde{I}_{1} + \tilde{I}_{2}\rightarrow\frac{1 - \e^{-\nu}}{1 -
\e^{2\beta- \nu}}.
\]

Similarly, one can write
\begin{eqnarray*}
\int_{0}^{{1}/{k_{n}}} Q(s) \,\mathrm{d}s &=& \frac{k_{n}}{A_{k_{n}}}\int
_{0}^{h_{2}(l(n)-1)\e^{-\nu(l(n)-1)}} Q(s) \,\mathrm{d}s - \frac
{k_{n}}{A_{k_{n}}}\int
^{h_{2}(l(n)-1)\e^{-\nu
(l(n)-1)}}_{h_{1}(n-1)\e^{-\nu(n-1)}} Q(s) \,\mathrm{d}s\\
& :=&\tilde{I}_{2}^{*}
- \tilde{I}_{1}^{*}.
\end{eqnarray*}
As shown in the case $0<\alpha<1$, we again use two different
representations of $l(n)$ for two different cases. Note that $\tilde
{I}_{1}^{*}$ is exactly $I_{1}^{*}$ considered in that case.

Observe that
\[
\tilde{I}_{2}^{*} =-\frac{k_{n}}{A_{k_{n}}} \sum
_{m= l(n)-1}^{\infty} L\bigl(\e^{2m+1}
\bigr)\e^{(2m+1)\beta}\bigl(h_{2}(m)\e^{-\nu m} -
h_{2}(m+1)\e^{-\nu
(m+1)}\bigr).
\]
As $n\rightarrow\infty$,
\[
\tilde{I}_{2}^{*}\rightarrow-c_{1}
\frac{(1- \e^{-\nu})\e^{\nu- \beta
}}{1 - \e^{2\beta- \nu}}\e^{(2\beta- \nu)\lceil({1}/{\nu})\log
c_{1}\rceil}
\]
and, from the case $0<\alpha<1$,
\[
\tilde{I}_{1}^{*}\rightarrow-\e^{\beta(2\lceil({1}/{\nu})\log
c_{1}\rceil- 1)}
\bigl(c_{1}\e^{-\nu\lceil({1}/{\nu})\log c_{1}\rceil} - 1\bigr).
\]

Finally, gathering the results above, we deduce the convergence to
$-\zeta$ where $\zeta$ is given by~(\ref{eeta1}).
\end{pf}

Theorem~\ref{tmain-th} concerns the partial sums $\sum_{j=1}^{n}
X_{j}$ along a subsequence $k_{n}$ of $n$. The following result
describes the behavior of the partial sums across all $n$. The result
is a direct consequence of Lemma $5$ of Meerschaert and Scheffler \cite
{meerschaertscheffler1998}. Recall that a collection of random
variables $\{Y_{n}\}_{n\geq1}$ is called \textit{stochastically compact}
if every subsequence $\{n^{\prime}\}$ has a further subsequence $\{n''\}
\subset\{n'\}$ for which $\{Y_{n''}\}$ converges in distribution. The
following notation will also be used. For a semi-stable distribution
$\tau$ with characteristic function $\psi(t)$, $\tau^{\lambda}$ will
denote the semi-stable distribution with the characteristic function
$\psi(t)^{\lambda}$.

\begin{proposition}\label{pmeerschaert-sheffler}
Let $X, X_{1}, X_{2},\ldots$ be i.i.d. random variables such that
%
\begin{equation}
\label{enewadded3} \frac{1}{A_{k_{n}}} \Biggl\{\sum_{j=1}^{k_{n}}X_{j}
- B_{k_{n}} \Biggr\} \mathop{\rightarrow}^{d} Y,
\end{equation}
where $Y$ follows a semi-stable distribution $\tau$ with $0<\alpha<2$
and $k_{n}$, $A_{k_{n}}$, $B_{k_{n}}$ are given in (\ref{econdkn}),
(\ref{econditionAkn}) and (\ref{edefiBkn}). Then, there exist
$a_{n}$ and $b_{n}$ such that $a_{n}$ is regularly varying with index
$\frac{1}{\alpha}$, $a_{k_{n}} = A_{k_{n}}$ and $a_{n}^{-1}(X_{1} +
X_{2}+\cdots+X_{n}) - b_{n}$ is stochastically compact, with every
limit point of the form $\lambda^{-{1}/{\alpha}}\tau^{\lambda}$ for
some $\lambda\in[1, c]$. Moreover, one can take
%
\begin{equation}
\label{enewadded4} a_{n} = \lambda_{n}^{{1}/{\alpha}}A_{k_{p_{n}}}
\quad\mbox {and}\quad b_{n} = \lambda_{n}^{1 - {1}/{\alpha}}
\frac{B_{k_{p_{n}}}}{A_{k_{p_{n}}}},
\end{equation}
where $\lambda_{n} = \frac{n}{k_{p_{n}}}$ and $p_{n}$, $k_{p_{n}}$ are
chosen so that $k_{p_{n}}\leq n < k_{p_{n+1}}$ for every $n \geq1$.
\end{proposition}

\begin{pf} The proposition follows directly from Lemma $5$ and its
proof in Meerschaert and Scheffler \cite{meerschaertscheffler1998}.
The left-hand side of (\ref{enewadded3}) appears in (2.9) of
Meerschaert and Scheffler \cite{meerschaertscheffler1998} as
\[
\tilde{a}_{n}^{-1}(X_{1}+X_{2}+
\cdots+X_{k_{n}}) - \tilde {b}_{n}.
\]
The existence of a regularly varying $a_{n}$ with $a_{k_{n}} = \tilde
{a}_{n}$ is part of the statement of Lemma $5$ of Meerschaert and
Scheffler \cite{meerschaertscheffler1998}. The expressions in (\ref
{enewadded4}) can be found in the proof of that Lemma $5$.
\end{pf}

\begin{corollary}\label{cmeerschaert-scheffler}
Under the assumptions of Proposition~\ref{pmeerschaert-sheffler},
%
\begin{equation}
\label{erighttailpro} \limsup_{n} P\bigl(a_{n}^{-1}(X_{1}+X_{2}+
\cdots+X_{n}) - b_{n} > x\bigr) \leq \sup
_{1\leq\lambda\leq c} P(Y_{\lambda} > x)
\end{equation}
and
%
\begin{equation}
\label{elefttailpro} \limsup_{n} P\bigl(a_{n}^{-1}(X_{1}+X_{2}+
\cdots+X_{n}) - b_{n} < x\bigr) \leq \sup
_{1\leq\lambda\leq c} P(Y_{\lambda} < x),
\end{equation}
where $Y_{\lambda}$ has the distribution of the form $\lambda^{-{1}/{\alpha}}\tau^{\lambda}$.
\end{corollary}

\begin{pf} Along a subsequence $\{n(k)\}$ of $\{n\}$, we have
\begin{eqnarray}
\label{eeq1}&& \limsup_{n} P\bigl(a_{n}^{-1}(X_{1}+X_{2}+
\cdots+X_{n}) - b_{n} > x\bigr)
\nonumber
\\[-8pt]
\\[-8pt]
\nonumber
&&\quad= \lim_{k}
P\bigl(a_{n(k)}^{-1}(X_{1}+X_{2}+
\cdots+X_{n(k)}) - b_{n(k)} > x\bigr).
\end{eqnarray}
Now, by Proposition~\ref{pmeerschaert-sheffler}, there exists a
further subsequence $\{n(k_{m})\}$ of $\{n(k)\}$ such that
%
\begin{eqnarray}
\label{eeq2} &&\lim_{m} P\bigl(a_{n(k_{m})}^{-1}(X_{1}+X_{2}+
\cdots+X_{n(k_{m})}) - b_{n(k_{m})} > x\bigr)
\nonumber
\\[-8pt]
\\[-8pt]
\nonumber
&&\qquad= P(Y_{\lambda} >
x),
\end{eqnarray}
where $Y_{\lambda}$ follows the distribution $\lambda^{-{1}/{\alpha
}}\tau^{\lambda}$. The relation (\ref{eeq2}) holds for all $x$ as long
as the semi-stable distribution $\tau^{\lambda}$ is continuous. By Huff
\cite{huff1972}, the continuity of $\tau^{\lambda}$ is equivalent to
$\int_{-\infty}^{0} \,\mathrm{d}L_{\lambda}(x) + \int_{0}^{\infty} \,\mathrm{d}R_{\lambda}(x)
= \infty$, where $L_{\lambda}$ and $R_{\lambda}$ define the L\'evy
measure of $\tau^{\lambda}$. By the definition of $\tau^{\lambda}$,
$L_{\lambda} = \lambda L$ and $R_{\lambda} = \lambda R$. Denote the
multiplicative period of $M_{L}(x)$ and $M_{R}(x)$ by $p>1$. Then,
after the change of variables $x = p^k y$ in the integrals below,
\begin{eqnarray*}
\int_{-\infty}^{0} \,\mathrm{d}L(x) + \int_{0}^{\infty}
\,\mathrm{d}R(x) &=& \sum_{k=-\infty
}^{\infty}\int
_{- p^{k+1}}^{-p^{k}}\,\mathrm{d}\frac{M_{L}(x)}{|x|^{\alpha}} + \sum
_{k=-\infty}^{\infty}\int_{p^{k}}^{p^{k+1}}
\,\mathrm{d}\frac
{(-M_{R}(x))}{x^{\alpha}}
\\
&=& \sum_{k=-\infty}^{\infty} p^{-k \alpha} \int
_{- p}^{-1} \,\mathrm{d}\frac
{M_{L}(y)}{|y|^{\alpha}} + \sum
_{k=-\infty}^{\infty} p^{-k \alpha} \int
_{1}^{p} \,\mathrm{d}\frac{(-M_{R}(y))}{y^{\alpha}} = \infty,
\end{eqnarray*}
unless $M_L \equiv0$ and $M_R \equiv0$.
Combining (\ref{eeq1}) and (\ref{eeq2}), we have (\ref
{erighttailpro}) for all $x\in\mathbb{R}$. The relation~(\ref
{elefttailpro}) can be obtained similarly.
\end{pf}

We will use Corollary~\ref{cmeerschaert-scheffler} to provide a
conservative confidence interval for $f_{W}(w)$ in Section~\ref{smain-results-sampling}.

\section{Application to sampling of finite point processes}
\label{smain-results-sampling}
We now turn back to the context of sampling of finite point processes.
The following result restates Theorem~\ref{tmain-th} and Proposition~\ref{pcentering} for the nonparametric estimator $\widehat{f}_{W}(w)$
of $f_{W}(w)$ given in
(\ref{efWestimator}) or (\ref{eempiricalfW})--(\ref{edescriptionxi}).

\begin{theorem}\label{tmain-th-sampling}
Suppose conditions (\ref{eevenWq})--(\ref{eh1ratio}) hold and $k_{n}$
is given in (\ref{edefinitionknAknBkn}). Let
%
\begin{equation}
\label{ealphalpha} \alpha= \frac{\nu}{2\log(q^{-1} - 1)}.
\end{equation}
If $\alpha\in(1,2)$, then
\[
d_{N}\bigl(\widehat{f}(w) - f(w)\bigr)\mathop{\rightarrow}^{d}
(-1)^{-w}(Y+\zeta),
\]
and
if $\alpha\in(0,1)$, then
\[
d_{N}\widehat{f}(w)\mathop{\rightarrow}^{d} (-1)^{-w}(Y+
\zeta),
\]
along the sample sizes $N = k_{n}$, where $d_{N} = \frac
{k_{n}}{A_{k_{n}}}$ with
%
\begin{equation}
\label{ecalculationAkn} A_{k_{n}}= \pmatrix{2n-2
\cr
w}(1-q)^{-w}
\bigl(q^{-1} - 1 \bigr)^{2n-2},
\end{equation}
and $\zeta$ defined in (\ref{eeta1}) and $Y$ is a semi-stable
distribution characterized by (\ref{echarcterization}) with
%
\begin{equation}
\label{ebeta-sampling} \beta= \log\bigl(q^{-1} - 1\bigr).
\end{equation}
\end{theorem}

\begin{pf}
In view of (\ref{eempiricalfW})--(\ref{edescriptionxi}), we are
interested in the distribution of
\[
X = \pmatrix{W_{q}
\cr
w}(-1)^{W_{q} - w} \frac{(1-q)^{W_{q} -
w}}{q^{W_{q}}}1_{\{W_{q} \geq w\}},
\]
where $w > 0$ is fixed and $W_{q}$ follows a p.m.f. satisfying (\ref
{eevenWq})--(\ref{eh1ratio}). For $W_{q} > w$ large enough, one can
write $(-1)^{w}X = L(\e^{W_{q}})\e^{\beta W_{q}} (-1)^{W_{q}}$ as given
in Theorem~\ref{tmain-th} with
%
\begin{equation}
L(x) =\pmatrix{\log x
\cr
w}(1 - q)^{-w} = (1-q)^{-w}
\frac{\prod_{i=0}^{w-1}(\log x - i)}{w!}
\end{equation}
and $\beta= \log\frac{1 - q}{q} = \log(q^{-1} - 1)$. Observe that
$L(x)$ is an ultimately increasing slowly varying function. Hence, when
$\alpha\in(1,2)$, by using (\ref{eempiricalfW})--(\ref{edescriptionxi}) and applying Theorem~\ref{tmain-th} and Proposition~\ref{pcentering},
\[
\frac{k_{n}}{A_{k_{n}}}\bigl(\widehat{f}_{W}(w) - f_{W}(w)
\bigr) = d_{N}\bigl(\widehat {f}_{W}(w) -
f_{W}(w)\bigr)
\]
converges to a semi-stable distribution $(-1)^{-w}(Y+\zeta)$ with
$\alpha$ in (\ref{ealphalpha}) and $A_{k_{n}}$ in (\ref{ecalculationAkn}).
When $\alpha\in(0,1)$,
\[
\frac{k_{n}}{A_{k_{n}}}\widehat{f}_{W}(w) = d_{N}\widehat
{f}_{W}(w)
\]
converges to a semi-stable distribution $(-1)^{-w}(Y+\zeta)$ with
$\alpha$ in (\ref{ealphalpha}) and $A_{k_{n}}$ in (\ref
{ecalculationAkn}).
\end{pf}

The next result provides a conservative confidence interval for $f(w)$
based on $\widehat f(w)$ when $1<\alpha<2$. The finite-sample
performance of the confidence interval and related issues are
considered in Chaudhuri and Pipiras \cite{chaudhuripipiras2013}.

\begin{proposition} \label{pconfidence} Under the assumptions and
notation of Theorem~\ref{tmain-th-sampling}, suppose $\alpha\in(1,2)$.
For $\gamma\in(0, 1)$, set
%
\begin{equation}
\mathcal{C} = \bigl[\widehat{f}_{W}(w) - \tilde{b}_{N}x_{1 - {\gamma
}/{2}},
\widehat{f}_{W}(w) - \tilde{b}_{N}x_{{\gamma}/{2}}\bigr],
\end{equation}
where
%
\begin{equation}
\tilde{b}_{N}= N^{{1}/{\alpha} - 1}A_{k_{p_{N}}}k_{p_{N}}^{-{1}/{\alpha}}
\end{equation}
with $p_{N}$ such that $k_{p_{N}}\leq N < k_{p_{N+1}}$ and
%
\begin{equation}
\sup_{1\leq\lambda\leq c}P\bigl(Y_{\lambda}^{\zeta} <
x_{{\gamma}/{2}}\bigr) = \frac{\gamma}{2},\qquad \sup_{1\leq\lambda\leq c}P
\bigl(Y_{\lambda}^{\zeta} > x_{1 - {\gamma
}/{2}}\bigr) =
\frac{\gamma}{2},
\end{equation}
where $Y_{\lambda}^{\zeta}$ has the distribution of the form $\lambda
^{-{1}/{\alpha}}\tau^{\lambda}$ and $\tau$ is the distribution of
$Y+\zeta$.
Then,
%
\begin{equation}
\label{econfi} \liminf_{N\rightarrow\infty} P\bigl(f_{W}(w) \in
\mathcal{C}\bigr) \geq1-\gamma,
\end{equation}
that is, $\mathcal{C}$ is a conservative $100(1-\gamma)\%$ confidence
interval for $f_{W}(w)$.
\end{proposition}
\begin{pf} When $\alpha\in(1,2)$, by using Corollary~\ref
{cmeerschaert-scheffler} and Theorem~\ref{tmain-th-sampling}, we get
\begin{eqnarray*}
&&\limsup_{N\rightarrow\infty}P\biggl(\frac{N\lambda_{N}^{-{1}/{\alpha
}}}{A_{k_{p_{N}}}}
\widehat{f}_{W}(w) - \lambda_{N}^{1 - {1}/{\alpha
}}
\frac{k_{p_{N}}}{A_{k_{p_{N}}}}f_{W}(w) < x_{{\gamma}/{2}}\biggr) \leq \sup
_{1\leq\lambda\leq c}P\bigl(Y_{\lambda}^{\eta} <
x_{{\gamma}/{2}}\bigr) = \frac{\gamma}{2}
\\
&&\Leftrightarrow\quad\limsup_{N\rightarrow\infty} P\biggl(\frac{N}{\lambda
_{N}k_{p_{N}}}
\widehat{f}_{W}(w) - \frac{\lambda_{N}^{{1}/{\alpha}
- 1}A_{k_{p_{N}}}}{k_{p_{N}}}x_{{\gamma}/{2}} <
f_{W}(w)\biggr) \leq\frac
{\gamma}{2}.
\end{eqnarray*}
Using $\lambda_{N}=\frac{N}{k_{p_{N}}}$, we get
%
\begin{equation}
\label{elefttailconfidence} \limsup_{N\rightarrow\infty}
P\bigl(\widehat{f}_{W}(w)
- N^{ {1}/{\alpha
} - 1}A_{k_{p_{N}}}k_{p_{N}}^{-{1}/{\alpha}}x_{
\gamma}/ {2} < f_{W}(w)\bigr)\leq\frac{\gamma}{2}.
\end{equation}
Similarly for the right tail, we get
%
\begin{eqnarray}
\label{erighttailconfidence} &&\limsup_{N\rightarrow\infty}
P\biggl(\frac{N\lambda_{N}^{-{1}/{\alpha
}}}{A_{k_{p_{N}}}}
\widehat{f}_{W}(w) - \lambda_{N}^{1 - {1}/{\alpha
}}
\frac{k_{p_{N}}}{A_{k_{p_{N}}}}f_{W}(w) > x_{1- {\gamma}/{2}}\biggr) \leq \sup
_{1\leq\lambda\leq c}P\bigl(Y_{\lambda}^{\eta} >
x_{1 -
{\gamma}/{2}}\bigr) = \frac{\gamma}{2}
\nonumber
\\
&&\quad\Leftrightarrow\quad\limsup_{N\rightarrow\infty}P\biggl(\frac{N}{\lambda
_{N}k_{p_{N}}}
\widehat{f}_{W}(w) - \frac{\lambda_{N}^{{1}/{\alpha}
- 1}A_{k_{p_{N}}}}{k_{p_{N}}}x_{1 - {\gamma}/{2}} >
f_{W}(w)\biggr) \leq \frac{\gamma}{2}
\\
&&\quad\Leftrightarrow\quad\limsup_{N\rightarrow\infty}P\bigl(\widehat{f}_{W}(w)
- N^{ {1}/{\alpha} - 1}A_{k_{p_{N}}}k_{p_{N}}^{-{1}/{\alpha
}}x_{1- {\gamma}/{2}}
> f_{W}(w)\bigr)\leq\frac{\gamma}{2}.\nonumber
\end{eqnarray}
Combining (\ref{elefttailconfidence}) and (\ref
{erighttailconfidence}), we get (\ref{econfi}).
\end{pf}

We conclude with two examples illustrating Theorem~\ref
{tmain-th-sampling}.

\begin{example}\label{exgeom}
Consider the case where $W$ follows a geometric distribution, that is,
$f_{W}(w) = c^{w-1}(1 - c)$, $w=1,2,3,\ldots$ and $0 < c < 1$.
Substituting this into (\ref{efwq-fw}) leads to
%
\begin{eqnarray}
\label{efwqgeometric} f_{W_{q}}(s) &=& \sum_{w = s}^{\infty}
\pmatrix{w
\cr
s}q^{s}(1-q)^{w-s}c^{w-1}(1-c).
\end{eqnarray}
When $s=0$, we get
%
\begin{equation}
\label{egeometric0} f_{W_{q}}(0) = \sum_{w=1}^{\infty}(1-q)^{w}c^{w-1}(1-c)
=\frac
{(1-q)(1-c)}{1 - c(1-q)}.
\end{equation}
When $s\geq1$, on the other hand, we have
%
\begin{eqnarray}
\label{egeometricsgeq1} f_{W_{q}}(s) &=& q^{s}c^{s-1}(1-c)\sum
_{w=s}^{\infty}\pmatrix{w
\cr
s}\bigl(c(1-q)
\bigr)^{w-s}
\nonumber
\\[-8pt]
\\[-8pt]
\nonumber
&=&\frac{q^{s}c^{s-1}(1-c)}{(1 - c(1-q))^{s+1}} = \frac
{c_{q}}{c}c_{q}^{s-1}(1
- c_{q}),
\end{eqnarray}
where $c_{q} = \frac{qc}{1 - c(1-q)}$, by using the identity $\sum_{w=s}^{\infty}{w\choose s} x^{w-s} = \sum_{r=0}^{\infty}{s+r\choose r}
x^{r} = (1 - x)^{-(s+1)}$.
Hence, for $x \geq1$,
%
\begin{equation}
\label{erighttailgeometric1} P \biggl(\frac{W_{q}}{2} \geq x, W_{q}\mbox{ is
even} \biggr) = \sum_{s=\lceil x \rceil}^{\infty}
\frac{c_{q}}{c}c_{q}^{2s-1}(1 - c_{q}) =
\frac{1}{c}\frac{c_{q}^{2\lceil x\rceil}}{1 + c_{q}}
\end{equation}
and
%
\begin{equation}
\label{elefttailgeometric1} P \biggl(\frac{W_{q}-1}{2} \geq x, W_{q}\mbox{ is
odd} \biggr) = \sum_{s=\lceil x\rceil}^{\infty}
\frac{c_{q}}{c}c_{q}^{2s}(1 - c_{q}) =
\frac
{c_{q}}{c}\frac{c_{q}^{2\lceil x \rceil}}{1+c_{q}}.
\end{equation}

Thus, the conditions (\ref{eevenWq})--(\ref{eh1ratio}) in Theorem~\ref
{tmain-th} are satisfied with $\nu= 2 \log\frac{1}{c_{q}}$,
$h_{1}(\lceil x \rceil) = \frac{1}{c(1+c_{q})}$, $h_{2}(\lceil x \rceil
) = \frac{c_{q}}{c(1+c_{q})}$ with $\frac{h_{2}(x)}{h_{1}(x)} = c_{q}$.
By using the expression of $\beta$ in (\ref{ebeta-sampling}), the
parameter $\alpha$ appearing in (\ref{echarcterization0}) or (\ref
{ealphalpha}) is given by
\[
\alpha= \frac{\log({1}/{c_{q}})}{\log(q^{-1} - 1)} = \frac{\log
{(1 - c(1-q))}/{(cq)}}{\log(q^{-1} - 1)}.
\]
Note that $c_{q} < 1$ and hence $\log\frac{1}{c_{q}} > 0$. Then,
$\alpha> 0$ is possible only when $q\in(0, 0.5)$. In particular, for
$q\in(0, 0.5)$,
%
\begin{eqnarray}
\label{eboundsgeometric1} 1 < \alpha< 2 \quad&\Leftrightarrow&\quad\frac{q}{1-q} < c <
\frac{1}{2(1-q)},
\\
\label{eboundsgeometric2} 0 < \alpha< 1\quad &\Leftrightarrow&\quad\frac{1}{2(1-q)} < c < 1.
\end{eqnarray}
Theorem~\ref{tmain-th-sampling} can now be applied in these two cases with
\[
A_{k_{n}} = {\pmatrix{2n-2
\cr
w}} (1-q)^{-w}
\bigl(q^{-1} - 1 \bigr)^{2n-2} \quad\mbox{and}\quad k_{n} =
\biggl\lceil\frac{c(1+c_{q})}{c_{q}^{2n - 2}} \biggr\rceil.
\]
\end{example}

\begin{remark*}
Under (\ref{eboundsgeometric1}) or (\ref
{eboundsgeometric2}), and $q\in(0, 0.5)$, the limit of $\widehat
{f}(w)$ involves a semi-stable distribution. On the other hand, as
proved in Antunes and Pipiras \cite{antunespipiras2011}, $\widehat
{f}(w)$ is asymptotically normal if $R_{q,w}<\infty$, where $R_{q,w}$
is given in (\ref{eRqw}). This condition obviously holds when $q\in
(0.5,1)$ (and also for $q = 0.5$ by recalling from Example~\ref
{exgeom} above that $f_{W_{q}}(s)\sim Cc_{q}^{s}\mbox{ as }s\rightarrow
\infty$). To understand when $R_{q,w}<\infty$ for $q\in(0,0.5)$,
observe that
%
\begin{eqnarray}
\label{egeometricvariance} R_{q,w}&=& \sum_{k=w}^{\infty}f_{W}(k)
(1-q)^{k-2w}\pmatrix{k
\cr
w}\sum_{s=w}^{k}
\pmatrix{s
\cr
w}\pmatrix{k-w
\cr
s-w} \biggl(\frac{1}{q} - 1
\biggr)^{s}
\nonumber
\\[-8pt]
\\[-8pt]
\nonumber
&=&\sum_{s=w}^{\infty}\pmatrix{s
\cr
w}
\bigl(q^{-1} - 1\bigr)^{s}\sum_{k=s}^{\infty
}c^{k-1}(1-c)
(1-q)^{k-2w}\pmatrix{k
\cr
w}\pmatrix{k-w
\cr
s-w}.
\end{eqnarray}
Since
\begin{eqnarray*}
\pmatrix{k
\cr
w}\pmatrix{k-w
\cr
s-w} =\frac{k!}{w!(k-w)!}\frac
{(k-w)!}{(s-w)!(k-s)!} =
\frac{k!}{(k-s)!s!}\frac{s!}{w!(s-w)!} = \pmatrix{k
\cr
s}\pmatrix{s
\cr
w},
\end{eqnarray*}
we have
%
\begin{eqnarray}
\label{egeometricvariancesecond} R_{q,w}&=&(1-c)\sum_{s=w}^{\infty}
\pmatrix{s
\cr
w}^{2}\bigl(q^{-1} - 1\bigr)^{s}\sum
_{k=s}^{\infty}\pmatrix{k
\cr
s}c^{k-1}(1-q)^{k-2w}
\nonumber
\\
&=&(1-c)\sum_{s=w}^{\infty}\pmatrix{s
\cr
w}^{2}\bigl(q^{-1} - 1\bigr)^{s}\sum
_{k=s}^{\infty}\pmatrix{k
\cr
s}\bigl(c(1-q)
\bigr)^{k-s}c^{s-1}(1-q)^{s-2w}
\nonumber
\\
&=& (1-c)\sum_{s=w}^{\infty}\pmatrix{s
\cr
w}^{2}\bigl(q^{-1} - 1\bigr)^{s}c^{s-1}(1-q)^{s-2w}
\sum_{k=s}^{\infty}\pmatrix{k
\cr
k-s}\bigl(c(1-q)
\bigr)^{k-s}
\nonumber
\\[-8pt]
\\[-8pt]
\nonumber
&=&(1-c)\sum_{s=w}^{\infty}\pmatrix{s
\cr
w}^{2}\bigl(q^{-1} - 1\bigr)^{s}c^{s-1}(1-q)^{s-2w}
\bigl(1 - c(1-q)\bigr)^{-(s+1)}
\\
&=& \biggl(\frac{1-c}{c} \biggr)\sum_{s=w}^{\infty}
\pmatrix{s
\cr
w}^{2}\bigl(q^{-1} - 1\bigr)^{s}
\bigl(c(1-q)\bigr)^{s}(1-q)^{-2w}\bigl(1 - c(1-q)
\bigr)^{-(s+1)}
\nonumber
\\
&=&d_{w}\sum_{s=w}^{\infty}\pmatrix{s
\cr
w}^{2} \biggl(\bigl(q^{-1} - 1\bigr)\frac
{c(1-q)}{1 - c(1-q)}
\biggr)^{s},\nonumber
\end{eqnarray}
where
$d_{w} =  (\frac{1-c}{c} )(1-q)^{-2w}\frac{1}{1 - c(1-q)}$.
Thus, $R_{q,w}<\infty$ if and only if
%
\begin{equation}
\label{egeometricvariancethird} \bigl(q^{-1} - 1\bigr)\frac{c(1-q)}{1 - c(1-q)} < 1
\quad\Leftrightarrow\quad c < \frac{q}{1-q}.
\end{equation}
Apart from the boundary cases $c = \frac{q}{1-q}$ and $c = \frac
{1}{2(1-q)}$, the ranges of $c$ given in (\ref{eboundsgeometric1}),
(\ref{eboundsgeometric2}) and (\ref{egeometricvariancethird}) now
cover the whole permissible interval $c\in(0,1)$.
\end{remark*}

\begin{example}\label{exneg-bin}
Consider the case where $W$ follows a negative binomial distribution,
that is, $f_{W}(w) = {w-1\choose r-1} c^{w-r}(1 - c)^{r}$, $w = r,
r+1,\ldots,0 < c < 1$. We first compute $f_{W_{q}}(s)$. One can
write $W = G_{1} + G_{2} +\cdots+ G_{r}$, where $G_{1}, G_{2},\ldots,
G_{r}$ are i.i.d. geometric random variables with p.m.f.
$f_{G_{1}}(w) = c^{w-1}(1-c), w\geq1$, and hence $W_{q} = G'_{1} +
G'_{2}+\cdots+G'_{r}$, where $G'_{1}, G'_{2},\ldots,G'_{r}$ are i.i.d.
random variables following the distribution given in (\ref
{egeometric0})--(\ref{egeometricsgeq1}).
Hence,
%
\begin{equation}
\label{enb0} f_{W_{q}}(0) = \biggl\{\frac{(1-q)(1-c)}{1 - c(1-q)} \biggr
\}^{r}.
\end{equation}
For $s\geq1$, we have
\[
f_{W_{q}}(s) = \sum_{i_{1}, i_{2},\ldots,i_{r}\geq0, i_{1} +
i_{2}+\cdots+i_{r} = s} P
\bigl(G^{\prime}_{1} = i_{1}\bigr)P
\bigl(G^{\prime}_{2} = i_{2}\bigr)\cdots P
\bigl(G^{\prime}_{r} = i_{r}\bigr).
\]
To evaluate this quantity, let
%
\begin{equation}
p_{j}^{r} = \sum_{i_{j+1}, i_{j+2},\ldots,i_{r}\geq1, i_{j+1} +
i_{j+2}+\cdots+i_{r} = s} P
\bigl(G^{\prime}_{j+1} = i_{j+1}\bigr)P
\bigl(G^{\prime}_{j+2} = i_{j+2}\bigr)\cdots P
\bigl(G^{\prime}_{r} = i_{r}\bigr),
\end{equation}
for $0\leq j < r$. Then, by using (\ref{egeometric0}),
\[
f_{W_{q}}(s) =\sum_{j=0}^{r-1}
\pmatrix{r
\cr
j} \biggl\{\frac{(1-q)(1-c)}{1 -
c(1-q)} \biggr\}^{j}
p_{j}^{r}.
\]
Now, by using (\ref{egeometricsgeq1}),
\begin{eqnarray*}
p_{j}^{r} &=& \biggl(\frac{c_{q}(1-c_{q})}{c}
\biggr)^{r-j}c_{q}^{s- (r-
j)}\sum
_{i_{j+1}, i_{j+2},\ldots,i_{r}\geq1, i_{j+1} + i_{j+2}+\cdots
+i_{r} = s} 1
\\
&=& \biggl(\frac{1-c_{q}}{c} \biggr)^{r-j}c_{q}^{s}
\pmatrix{s-1
\cr
r-j-1}.
\end{eqnarray*}
Hence, for $s\geq1$,
\begin{eqnarray*}
f_{W_{q}}(s) &=& c_{q}^{s}\sum
_{j=0}^{r-1} \biggl\{\frac{(1-q)(1-c)}{1 -
c(1-q)} \biggr
\}^{j} \biggl(\frac{1-c_{q}}{c} \biggr)^{r-j}\pmatrix{r
\cr
j}\pmatrix{s-1
\cr
r-j-1}
\\
&=& c_{q}^{s-1}p^{*}(s),
\end{eqnarray*}
where $p^{*}(s)$ is a polynomial given as
\[
p^{*}(s) = \sum_{i=1}^{r - 1}a^{*}_{i}s^{i}.
\]

This implies that for $x > 1$,
%
\begin{equation}
\label{erighttailnb1} P \biggl(\frac{W_{q}}{2} \geq x, W_{q}\mbox{ is
even} \biggr) = \sum_{s=\lceil x \rceil}^{\infty}c_{q}^{2s-1}p^{*}(2s)
= c_{q}^{2\lceil x
\rceil}\sum_{s=\lceil x \rceil}^{\infty}
c_{q}^{2s - 2\lceil x \rceil
-1}p^{*}(2s)
\end{equation}
and
%
\begin{equation}
\label{elefttailnb1} P \biggl(\frac{W_{q} - 1}{2} \geq x, W_{q}\mbox{ is
odd} \biggr) = \sum_{s=\lceil x \rceil}^{\infty}c_{q}^{2s}p^{*}(2s+1)
= c_{q}^{2\lceil x
\rceil}\sum_{s= \lceil x \rceil}^{\infty}
c_{q}^{2s - 2\lceil x \rceil
}p^{*}(2s+1).
\end{equation}
Thus the conditions (\ref{eevenWq})--(\ref{eoddWq}) in Theorem~\ref
{tmain-th} are satisfied with $\nu= 2 \log\frac{1}{c_{q}}$,
$h_{1}(x) = \sum_{k=0}^{\infty}c_{q}^{2k - 1}p^{*}(2x+2k)$, $h_{2}(x) =
\sum_{k=0}^{\infty}c_{q}^{2k}p^{*}(2x+1+2k)$.
The conditions (\ref{eh2h1ratio})--(\ref{eh1ratio}) also hold with
$c_{1} = c_{q}$. The parameter $\alpha$ appearing in (\ref
{echarcterization0}) is given by
\[
\alpha= \frac{\log({1}/{c_{q}})}{\log(q^{-1} - 1)} = \frac{\log
({(1 - c(1-q))}/{(cq)})}{\log(q^{-1} - 1)}.
\]
Note that $c_{q} < 1$ and hence $\log\frac{1}{c_{q}} > 0$. Then,
$\alpha> 0$ is possible only when $q\in(0, 0.5)$. In particular, for
$q\in(0, 0.5)$, the two cases (\ref{eboundsgeometric1})--(\ref
{eboundsgeometric2}) can be considered.
Theorem~\ref{tmain-th-sampling} can now be applied in these two cases with
\[
A_{k_{n}} = {\pmatrix{2n-2
\cr
w}}(1-q)^{-w}
\bigl(q^{-1} - 1 \bigr)^{2n-2} \quad\mbox{and}\quad k_{n} =
\biggl\lceil\frac{1}{c_{q}^{2n - 2}h_{1}(n-1)} \biggr\rceil.
\]
\end{example}

\begin{appendix}\label{sauxiliary}

\section*{Appendix: Auxiliary results}

We state and prove here a number of auxiliary results used in Section~\ref{smain-results}.

\begin{lemma}\label{lg1-tilde-g1}
Let $g_{1}$ and $\tilde{g}_{1}$ be defined in (\ref{edefinitiong1})
and (\ref{edefinitiong1tilde}), respectively. Then, $\tilde{g}_{1}(y)
- g_{1}(y)\rightarrow0, \mbox{ as } y\rightarrow\infty$.
\end{lemma}

\begin{pf} For $n\geq2$, if
\[
\label{eintervaltildeg1g1differcase1} n - 1 + \frac{1}{2\beta}\log L\bigl(\e^{2n-2}\bigr)\leq
y < n - 1 + \frac
{1}{2\beta}\log L \bigl(\e^{2n} \bigr),
\]
then
\setcounter{equation}{0}
\begin{equation}
\label{etildeg1g1differcase1} 0 \leq\tilde{g}_{1}(y) - g_{1}(y) <
\frac{1}{2\beta}\log \frac{L
(\e^{2n} )}{L (\e^{2n-2} )}\rightarrow0\qquad\mbox{as } y\rightarrow
\infty\ (n\rightarrow\infty),
\end{equation}
since $L$ is a slowly varying function.
If
\[
\label{eintervaltildeg1g1differcase2} n - 1 + \frac{1}{2\beta}\log L\bigl(\e^{2n}\bigr)\leq
y < n + \frac{1}{2\beta
}\log L \bigl(\e^{2n} \bigr),
\]
then similarly
%
\begin{equation}
\label{etildeg1g1differcase2} \tilde{g}_{1}(y) - g_{1}(y) =
\frac{1}{2\beta}\log\frac{L
(\e^{2n} )}{L (\e^{2n-2} )}\rightarrow0\qquad\mbox{as } y\rightarrow
\infty\ (n\rightarrow\infty).
\end{equation}
\end{pf}

\begin{lemma}\label{lg-tilde-*}
Let $\tilde{g}_{1}^{*}$ be defined in (\ref{etildeg1*}). Then, for any
$A > 0$,
\[
\label{esecondlemma} \tilde{g}_{1}^{*} (\log Ax ) -
\tilde{g}_{1}^{*} (\log x )\rightarrow0 \qquad\mbox{as } x
\rightarrow\infty.
\]
\end{lemma}
\begin{pf}Suppose without loss of generality that $A > 1$. First,
note that
%
\begin{eqnarray}
\label{eg1tildestarproperty} \tilde{g}_{1}^{*} (\log Ax ) -
\tilde{g}_{1}^{*} (\log x ) &=& \frac{1}{2\beta} \bigl(
\log L \bigl(\e^{2n_{Ax} - 2} \bigr) - \log L \bigl(\e^{2n_{x} - 2} \bigr)
\bigr)
\nonumber
\\
&=& \frac{1}{2\beta}\log\frac{L (\e^{n_{Ax} - 2} )}{L
(\e^{n_{x} - 2} )}
\\
&=& \frac{1}{2\beta}\log\frac{L (\e^{2n_{Ax} - 2n_{x}} \e^{2n_{x} -
2} )}{L (\e^{n_{x} - 2} )},\nonumber
\end{eqnarray}
where, for $y$ ($=x$ or $Ax$),
\[
n_{y}-1 + \frac{1}{2\beta}\log L\bigl(\e^{2n_{y}-2}\bigr) \leq
\log y < n_{y} + \frac{1}{2\beta}\log L\bigl(\e^{2n_{y}}
\bigr).
\]
Observe that $n_{Ax} - n_{x}$ takes only positive integer values, and that
\[
0 \leq n_{Ax} - n_{x} \leq\lceil\log A\rceil.
\]
Hence, by Theorem~1.2.1 of Bingham, Goldie and Teugels \cite
{binghamgoldieteugels1987},
\begin{eqnarray*}
\frac{L (\e^{n_{Ax} - n_{x}} \e^{n_{x} - 1} )}{L (\e^{n_{x} -
1} )}\rightarrow1 \qquad\mbox{as } \e^{n_{x} - 1} \rightarrow\infty
\ (\mbox{or } x\rightarrow\infty).
\end{eqnarray*}
This yields the result.
\end{pf}

\begin{lemma}\label{ll1*}
The function $l_{1}^{*}(x)$ defined in (\ref{el1*x}) is
right-continuous and slowly varying at $\infty$.
\end{lemma}

\begin{pf} To show that $l_{1}^{*}(x)$ is slowly varying, write
\begin{eqnarray*}
l_{1}^{*}(x) &=&\frac{h_{1} (\lceil g_{2}(({1}/{(2\beta)})\log
x)\rceil_{\inplus} )}{h_{1} (g_{2}(({1}/{(2\beta)})\log x)
)}\\
&&{}\times h_{1}
\biggl(g_{2}\biggl(\frac{1}{2\beta}\log x\biggr)
\biggr)\e^{\nu\tilde
{g}_{1}^{*}(({1}/{(2\beta)})\log x)}\e^{-\nu({g}_{1}(({1}/{(2\beta)
})\log x) - \tilde{g}_{1}(({1}/{(2\beta)})\log x))}.
\end{eqnarray*}
Note that
\begin{eqnarray*}
\label{eratioinlemma} &&\frac{h_{1} (\lceil g_{2}(({1}/{(2\beta)})\log x)\rceil
_{\inplus} )}{h_{1} (g_{2}(({1}/{(2\beta)})\log x) )}\\
&&\quad = \frac{h_{1} (({\lceil g_{2}(({1}/{(2\beta)})\log x)\rceil
_{\inplus}}/{(g_{2}(({1}/{(2\beta)})\log x))})g_{2}(({1}/{(2\beta)})\log
x) )}{h_{1} (g_{2}(({1}/{(2\beta)})\log x) )}\rightarrow 1
\end{eqnarray*}
by using (\ref{eh1ratio}), since $g_{2}(\frac{1}{2\beta}\log
x)\rightarrow\infty$ and
\[
\frac{\lceil g_{2}(({1}/{(2\beta)})\log x)\rceil_{\inplus}}{g_{2}(
({1}/{(2\beta)})\log x)}\rightarrow1 \qquad\mbox{as } x\rightarrow\infty .
\]
By Lemma~\ref{lg1-tilde-g1}, we also have
\[
\label{esecondpartlemma3} \e^{-\nu({g}_{1}(({1}/{(2\beta)})\log x) -
\tilde{g}_{1}(({1}/{(2\beta)})\log x))}\rightarrow1\qquad \mbox{as }x\rightarrow\infty.
\]
Hence, $l_{1}^{*}(x)$ is asymptotically equivalent to
%
\begin{equation}
\label{easymptoticequivalent} h_{1} \biggl(g_{2}\biggl(\frac{1}{2\beta}
\log x\biggr) \biggr)\e^{\nu\tilde
{g}_{1}^{*}(({1}/{(2\beta)})\log x)}.
\end{equation}

It is enough to show that the function (\ref
{easymptoticequivalent}) is slowly varying. By using Lemma~\ref
{lg-tilde-*}, we have
%
\begin{equation}
\label{ethirdpartlemma3} \frac{\e^{\nu\tilde{g}_{1}^{*}(({1}/{(2\beta)})\log Ax)}}{\e^{\nu\tilde
{g}_{1}^{*}(({1}/{(2\beta)})\log x)}}\rightarrow1 \qquad\mbox{as }x\rightarrow\infty.
\end{equation}
It remains to show that $h_{1}(g_{2}(\frac{1}{2\beta}\operatorname{log }x))$ is
a slowly varying function. For $A > 0$,
%
\begin{eqnarray}
\label{elastpartlemma3} &&\frac{h_{1}(g_{2}(({1}/{(2\beta)})\log Ax))}{h_{1}(g_{2}
(({1}/{(2\beta))
}\log x))}
\nonumber
\\[-8pt]
\\[-8pt]
\nonumber
&&\quad= \frac{h_{1} (({g_{2}(({1}/{(2\beta)})\log
Ax)}/{(g_{2}(({1}/{(2\beta)})\log x)}))g_{2}(({1}/{(2\beta)})\log x)
)}{h_{1}(g_{2}(({1}/{(2\beta)})\log x))}.
\end{eqnarray}
Now, by using (\ref{eg2representation}),
\begin{eqnarray*}
\label{elastpartlemma3extra}&& \frac{g_{2}(({1}/{(2\beta)})\log Ax)}{g_{2}
(({1}/{(2\beta)})\log x)} \\
&&\quad= \frac{({1}/{(2\beta)})\log Ax + g_{2}^{*}(({1}/{(2\beta)})\log
Ax)}{({1}/{(2\beta)})\log x + g_{2}^{*}(({1}/{(2\beta)})\log
x)}
\\
&&\quad=1 + \frac{({1}/{(2\beta)})\log Ax + g_{2}^{*}(({1}/{(2\beta)})\log
Ax) - ({1}/{(2\beta)})\log x - g_{2}^{*}(({1}/{(2\beta)})\log x)}{
({1}/{(2\beta)})\log x + g_{2}^{*}(({1}/{(2\beta)})\log x)}
\\
&&\quad=1 + \frac{({1}/{(2\beta)})\log A + g_{2}^{*}(({1}/{(2\beta)})\log
Ax) - g_{2}^{*}(({1}/{(2\beta)})\log x)}{ g_{2}(({1}/{(2\beta)})\log
x)}\rightarrow1,
\end{eqnarray*}
since $g_{2}(\frac{1}{2\beta}\log x)\rightarrow\infty$ and by using
(\ref{eg2*property}), $g_{2}^{*}(\frac{1}{2\beta}\log Ax) -
g_{2}^{*}(\frac{1}{2\beta}\log x)\rightarrow0$. Thus, by using~(\ref
{eh1ratio}) and (\ref{elastpartlemma3}), we have
\begin{eqnarray*}
\label{elastpartlemma3proved} \frac{h_{1}(g_{2}(({1}/{(2\beta)})\log Ax))}{h_{1}(g_{2}(\frac{1}{2\beta
}\log x))}\rightarrow1\qquad\operatorname{as }x\rightarrow
\infty.
\end{eqnarray*}
This completes the proof that $l_{1}^{*}(x)$ is a slowly varying function.

The function $l_{1}^{*}(x)$ is right-continuous since $h_{1}(x)$
can be defined to be continuous, $g_{2}$ is continuous (as the inverse
of a continuous increasing function) and $g_{1}$, $\tilde{g}_{1}$ and
$\tilde{g}_{1}^{*}$ are right-continuous functions.
\end{pf}

\begin{lemma}\label{lfinite-number}
Let $L$ be a slowly varying function. Then, for any fixed $x_{0}\neq
\e^{2\beta(r+1 - b_{1})}$, $r\in\mathbb{Z}$, $\beta> 0$, there are only
finitely many integer values of $n$ for which
%
\begin{equation}
\label{econditionfinite} m - b_{1} +\frac{1}{2\beta}\log L
\bigl(\e^{2m - b_{2}}\bigr) \leq\frac{1}{2\beta
}\log(A_{k_{n}}x_{0})
< m - b_{1} +\frac{1}{2\beta}\log L\bigl(\e^{2m - b_{3}}\bigr),
\end{equation}
where $A_{k_{n}} =\e^{(n-1)2\beta}L(\e^{2n-2})$, $m$ takes positive
integer values, $b_{1}$, $b_{2}$ and $b_{3}$ are fixed positive
constants with $b_{2} > b_{3}$.
\end{lemma}

\begin{pf}
Suppose $m = n + r_{n}$, where $r_{n}$ is a sequence of integers. We
first show that if (\ref{econditionfinite}) is satisfied for
infinitely many values of $n$, then $\sup_{n\geq1}|r_{n}|<\infty$.
Arguing by contradiction, for example, assume $r_{n}\rightarrow\infty$
as $n\rightarrow\infty$. From (\ref{econditionfinite}), we need to have
%
\begin{equation}
\label{efirsteqlemma4} \e^{2\beta(r_{n} + 1 - b_{1})}\frac{L(\e^{2n + 2r_{n} -
b_{2}})}{L(\e^{2n-2})}\leq x_{0} <
\e^{2\beta(r_{n} + 1 - b_{1})}\frac
{L(\e^{2n + 2r_{n} - b_{3}})}{L(\e^{2n - 2})}.
\end{equation}
A standard argument using Potter's bounds for $L$ shows that $\e^{2\beta
(r_{n} + 1 - b_{1})}\frac{L(\e^{2n + 2r_{n} - b})}{L(\e^{2n -
2})}\rightarrow\infty$ $(b= b_{2} \mbox{ or } b_{3})$ when
$r_{n}\rightarrow\infty$. Since $x_{0}$ is fixed, this leads to a
contradiction. A similar argument can be applied when $r_{n}\rightarrow
-\infty$.

Next we show that $m$ is necessarily of the form $m = n+ r$ where
$r$ is a fixed integer for large enough $n$. We prove this by
contradiction. First, observe that $r_{n}$ can only take finitely many
integer values. Now if $r_{n}$ has a subsequence $r_{n_{k}}\rightarrow
r$, then letting $n\rightarrow\infty$ in (\ref{efirsteqlemma4}), we
have $\e^{2\beta(r+1 - b_{1})} = x_{0}$. Thus, $r$ is determined
uniquely and since $r_{n}$ are integers, we have that $r_{n} = r$ for
large enough $n$.

Finally, if $m = n+ r$, then (\ref{econditionfinite}) cannot hold
for infinitely many values of $n$ unless $x_{0} = \e^{2\beta(r + 1 -
b_{1})}$. This proves the lemma.
\end{pf}

%
\begin{lemma}\label{ltwo-slow-var}
Let (\ref{ecorollary3,megyesi,lefttail})--(\ref{ecorollary3,megyesi,righttail}) hold for a random variable $X$ with $l^{*}(x)$
replaced by a right-continuous slowly varying function $l_{1}^{*}(x)$
in (\ref{ecorollary3,megyesi,lefttail}). Then, $l^{*}(x)$ in (\ref
{ecorollary3,megyesi,righttail}) can be replaced by another
right-continuous function $l_{2}^{*}(x)$ if $\frac
{l_{2}^{*}(x)}{l_{1}^{*}(x)}\rightarrow1$ as $x\rightarrow\infty$.
\end{lemma}
\begin{pf}
Observe that
%
\begin{eqnarray}
\label{eupdatedone} 1 - F(x) &=& x^{-\alpha}l_{2}^{*}(x)
\bigl(M_{R}\bigl(\delta(x)\bigr) + h_{R}(x)\bigr)
\nonumber
\\
&=&x^{-\alpha}l_{1}^{*}(x) \biggl(M_{R}
\bigl(\delta(x)\bigr) + h_{R}(x) + \biggl(\frac
{l_{2}^{*}(x)}{l_{1}^{*}(x)} - 1
\biggr) \bigl(M_{R}\bigl(\delta(x)\bigr) + h_{R}(x)\bigr)
\biggr)
\\
&=&x^{-\alpha}l_{1}^{*}(x) \bigl(M_{R}
\bigl(\delta(x)\bigr) + h_{R}(x) + \tilde{h}_{R}(x)\bigr),\nonumber
\end{eqnarray}
where
%
\begin{equation}
\tilde{h}_{R}(x) = \biggl(\frac{l_{2}^{*}(x)}{l_{1}^{*}(x)} - 1\biggr)
\bigl(M_{R}\bigl(\delta (x)\bigr) + h_{R}(x)\bigr).
\end{equation}
Since $\frac{l_{2}^{*}(x)}{l_{1}^{*}(x)}\rightarrow1$ as $x\rightarrow
\infty$, $M_{R}$ is a bounded periodic function from (\ref{elevyleft})
and $h_{R}(A_{k_{n}}x)\rightarrow0$, as $n\rightarrow\infty$, we have
$\tilde{h}_{R}(A_{k_{n}}x)\rightarrow0$ for every continuity point $x$
of $M_{R}(x)$. Hence, in (\ref{eupdatedone}), one can take the new
error function to be $h_{R}(x) + \tilde{h}_{R}(x)$. Hence, the result
is proved.
\end{pf}

%
\begin{lemma}\label{lh1-bounds} Let $h_{1}$ be the function defined
in Theorem~\ref{tmain-th} and satisfying (\ref{eh1ratio}). For every
$\delta> 0$, there is $M_{\delta}$ such that, for all $n >M_{\delta}$,
\[
h_{1}(M_{\delta}+1)\e^{M_{\delta} + 1}\e^{\delta n}<h_{1}(n)<
\frac
{h_{1}(M_{\delta}+1)}{\e^{\delta(M_{\delta} + 1)}}\e^{\delta n}.
\]
\end{lemma}
\begin{pf}
Fix any $\delta=\delta_{0}\in(0,1)$. By using (\ref{eh1ratio}), there
exists $M_{\delta_{0}}$ such that for all $m > M_{\delta_{0}}$, $1 -
\delta_{0}<\frac{h_{1}(m+1)}{h_{1}(m)} < 1 + \delta_{0}$. Take any $n >
M_{\delta_{0}}$. Then,
\begin{eqnarray*}
h_{1}(n) &=& \frac{h_{1}(n)}{h_{1}(n-1)}\frac
{h_{1}(n-1)}{h_{1}(n-2)}\cdots
\frac{h_{1}(M_{\delta
_{0}}+2)}{h_{1}(M_{\delta_{0}}+1)}h_{1}(M_{\delta_{0}}+1)
\\
&<& h_{1}(M_{\delta_{0}}+1) (1+\delta_{0})^{n - M_{\delta_{0}} - 1}
< h_{1}(M_{\delta_{0}}+1)\e^{\delta_{0}(n - M_{\delta_{0}} - 1)}.
\end{eqnarray*}
Similarly,
\[
h_{1}(n) > h_{1}(M_{\delta_{0}}+1) (1 -
\delta_{0})^{n - M_{\delta_{0}}
- 1} > h_{1}(M_{\delta_{0}}+1)\e^{-\delta_{0}(n - M_{\delta_{0}} -
1)}.
\]
\upqed\end{pf}
\end{appendix}

\section*{Acknowledgments}

The authors would like to thank the Associate Editor and an
anonymous referee for
useful comments and suggestions.



%




\printhistory
\end{document}